\documentclass[12pt,draft]{amsart}

\usepackage{amssymb}

\usepackage{amsmath}
\usepackage{amsthm}
\usepackage{amscd}

\newtheorem{theorem}{Theorem}
\newtheorem{lem}{Lemma}
\newtheorem{prop}{Proposition}
\newtheorem{cor}{Corollary}

\theoremstyle{definition}

\newtheorem{rem}{Remark}

\newtheorem{exa}{Example}

\numberwithin{equation}{section} \numberwithin{rem}{section}
\numberwithin{exa}{section} \numberwithin{theorem}{section}
\numberwithin{lem}{section} \numberwithin{prop}{section}
\numberwithin{cor}{section} \numberwithin{Def}{section}

\def\1{1\!\!\hbox{{\rm I}}}
\renewcommand{\Re}{{\Bbb R}}
\def\ax{\Re^+}
\def\CC{{\Bbb C}}
\def\TT{{\Bbb T}}
\def\eps{\varepsilon}

\newcommand{\prt}{\partial}
\newcommand{\kap}{\varkappa}
\newcommand{\Af}{\mathcal{A}}
\newcommand{\Mf}{\mathcal{M}}
\newcommand{\Df}{\mathcal{D}}
\newcommand{\Kf}{\mathcal{K}}
\newcommand{\be}{\begin{equation}}
\newcommand{\ee}{\end{equation}}
\newcommand{\im}{\mathop{\mathrm{Im}}}
\newcommand{\re}{\mathop{\mathrm{Re}}}
\newcommand{\ba}{\begin{aligned}}
\newcommand{\ea}{\end{aligned}}

\begin{document}

\title[{Exact asymptotic for distribution
densities of L\'evy functionals}]
    {{Exact asymptotic for distribution
densities of L\'evy functionals}}

\author{Victoria P. Knopova$\,^*$}
\author{Alexey M. Kulik}
\address{Kiev 03187 Acad. Glushkov Ave. 40, V.M.Glushkov
Institute of Cybernetics, Ukrai\-ni\-an National Academy of
Sciences} \email{vic\_knopova@gmx.de}
\address{Kiev 01601 Tereshchenkivska str. 3, Institute of Mathematics,
Ukrai\-ni\-an National Academy of Sciences}
 \email{kulik@imath.kiev.ua}

\thanks{$\,^*$The  DAAD scholarship during June -- August 2009 is gratefully acknowledged.}

%    General info
%\date{02/09/2003}
%\dedicatory{We dedicate the paper to our teachers.}

\begin{abstract}
A version of the saddle point method is developed,  which allows {one} to describe exactly the asymptotic behavior of distribution densities of  L\'evy driven stochastic integrals with deterministic kernels. Exact asymptotic behavior is established for  (a) the transition
probability density of a real-valued L\'evy process; (b) the
transition probability density and the invariant distribution density of a L\'evy driven
Ornstein-Uhlenbeck process; (c) the distribution density of the fractional  L\'evy motion.
\end{abstract}

\maketitle

\section{Introduction}

In {this} paper, we develop a version of the {\it saddle point method},  which allows {one} to describe exactly the asymptotic behavior of distribution densities of  L\'evy processes and, more generally, L\'evy driven stochastic integrals with deterministic kernels.  We start the exposition with the outline of the principal idea of the approach.

Let $(Z_t)_{t\geq 0}$  be a real-valued L\'evy process  with
characteristic exponent $\psi$; that is,
\begin{equation}
E e^{i z Z_t} =e^{t\psi(z)}, \quad t>0.\label{char}
\end{equation}
The function $\psi:\mathbb{R}\to\mathbb{C}$ (the \emph{characteristic exponent} of the process $Z$) admits the
L\'evy-Khinchin representation
\begin{equation}\label{LKh}
\psi(z)=iaz - b z^2+\int_{\mathbb{R}}
 \left( e^{iuz}-1- iz u 1_{\{|u|\leq 1\}} \right)\mu(du),
\end{equation}
where $a\in\mathbb{R}$, $b\geq 0,$ and $\mu(\cdot)$ is a \emph{L\'evy
measure}, i.e. $\int_{\mathbb{R}} (1\wedge u^2)\mu(du)<\infty$.
 Under  some   conditions (see Section \ref{s2} below), the function
$e^{t\psi}$ is  integrable, and hence the transition probability
density  $p_t(x)$ of  the process $Z_t$ has the integral
representation as  the inverse Fourier transform of the
characteristic function \eqref{char}:
\begin{equation}
p_t(x)=\frac{1}{2\pi} \int_{\mathbb{R}} e^{-iz  x+t\psi(z)}\,
dz. \label{ptx}
\end{equation}
Our intent is to investigate the  oscillatory
integral  (\ref{ptx}) using  the {\it
saddle point method}. According
 to this method (see \cite{Co65}), one can,
 under the assumption that the characteristic exponent $\psi$ admits  an analytic
extension to the complex plane,  apply the Cauchy theorem in order
to change the  integration path  in \eqref{ptx}:
    \begin{equation}
    p_t(x)=\frac{1}{2\pi} \int_{C} e^{-iz  x+t\psi(z)}\, dz.\label{ptx1}
    \end{equation}
Here $C$ is  certain   properly chosen contour that  allows {one} to  apply
the \textit{Laplace method} (\cite{Co65}, \cite{Ev62},
\cite{Fe77}) for estimating integral (\ref{ptx1}). A perfect
choice of the contour $C$ would be  the proper branch of the
curve $\{z: \im(-iz x+t\psi(z))=\im(-iz_0 x + t\psi(z_0))\}$,
 where
$z_0$ is a critical point  of   the function $-izx+t\psi(z)$ ({\it a
saddle point}). Under such a choice the integrand in (\ref{ptx1})
is real-valued; in this case the {\it saddle point method}
coincides with  the {\it fastest descent method}, see
\cite{Co65}.  However the complicated ``oscillatory'' structure of
the L\'evy-Khinchin representation of  $\psi$ does not give  an
opportunity to solve  the equation  $\im(-iz x+t\psi(z))=\im(-iz_0
x+t\psi(z_0))$ explicitly.    Instead, we put in \eqref{ptx1}
$C=\Re+i\xi_0$ with $i\xi_0$ being a critical point  of  the
function $-izx+t\psi(z)$. Under such a choice, we develop an
appropriate version of the Laplace method and give exact
asymptotics for the transition probability density $p_t(x)$.

The saddle point method is a classic tool for
estimating  a distribution density in various versions of the
local limit theorem with  the normal domain of attraction (see
\cite{IL65}, chapters 8, 10, and the references therein). In the L\'evy
processes setting, the idea of applying the complex analysis
technique was used, for instance, in \cite{KS} for getting upper estimates for
\eqref{ptx} in  the  case  when the characteristic exponent is real valued.

Since we require the characteristic exponent $\psi$ to have an analytic
extension to the complex plane, a standing assumption on the L\'evy measure within our approach is that it is exponentially integrable; that is,
\begin{equation}
\int_{|y|\geq 1} e^{C  y} \mu(dy)<\infty\quad \text{for all
$C\in\mathbb{R}$}.
\label{esint}
\end{equation}
  Equivalently, \eqref{esint} means that the variable $Z_1$
has exponential moments, i.e. $Ee^{cZ_1}<\infty$ for all
$c\in\mathbb{R}$, see \cite{Sa99}, $\S25-26$. Assumption
\eqref{esint} is non-restrictive, and is satisfied, for instance,
for a {\itshape generalized   tempered} L\'evy measure  of the
form $\mu(du)=\psi(u)\tilde \mu(du)$, where $\tilde \mu$ is
another L\'evy measure, and $\psi$ has a super-exponential decay,
i.e., $e^{Cu}\psi(u)\to 0$, $u\to \infty$, for all
$C\in\mathbb{R}$. For various results on generalized  tempered
L\'evy processes and  models that lead to processes of such a
type, we refer {the} reader to  Rosinski and Singlair \cite{RS10},
Sztonyk \cite{St10a}, \cite{St10b},  Bianchi et. al.
\cite{BRKF08}, \cite{KRBF09}. The notion of a generalized
tempered L\'evy measure is closely related to the notions of  a
\textit{tempered} and a \textit{layered} Levy measure, with the
function $\psi(u)$ in the above definition  respectively being
completely monotonous or having a polynomial decay rate (for both
these classes (\ref{esint}) fails). For the results on tempered
and layered L\'evy  processes and related models, see Rosinski
\cite{Ro07}, Cohen and Rosinski \cite{CR07}, Cont and Tankov
\cite{CT04}, Carr et. al. \cite{CGMY02}, \cite{CGMY03}, Baeumer
and Meerschaert \cite{BM09}, Kim et. al. \cite{KRCB09}, Houdr\'e
and Kawai \cite{HK06}. Of course,  this list of references is far
from  complete.

The method described above can be extended naturally for L\'evy driven stochastic integrals with deterministic kernels. Let
\be\label{y}
Y_t:=\int_I f(t,s) dZ_s,
\ee
where $I\subset \Re$ is an interval, $f$ is a deterministic function, and $Z_t$ is a L\'evy process
(in some particularly important cases, one  should take   $I=\Re$,    and then $Z$ should be assumed to
be two-sided; see details in Section \ref{s2} below).  The characteristic exponent of $Y_t$ can be
written explicitly (see \eqref{chint} below),   which makes it possible to apply the  method described
above to study the asymptotic behaviour of the distribution density of $Y_t$.

We mention two particular classes of processes,    frequently used in applications,
and having representation (\ref{y}).  The \emph{L\'evy driven Ornstein-Uhlenbeck process} is
defined as the solution to the linear SDE
\begin{equation}
dX_t =\gamma X_t\, dt +d Z_t,\quad t\geq 0,\label{sde}
\end{equation}
and has the integral representation
\be\label{sde1}
X_t=e^{\gamma t}X_0+\int_0^te^{\gamma(t-s)}\,dZ_s, \quad t\geq 0.
\ee
 If the initial value $X_0$ is non-random,   the distributional properties of $X_t$ are determined
 by the second term in the right hand side of (\ref{sde1}), which clearly has the form (\ref{y})
 with $I=\ax$ and $f(t,s)=e^{\gamma(t-s)}\1_{s\leq t}$. In what follows, we call such a process a
 \emph{non-stationary version of the Ornstein-Uhlenbeck process}.

The Ornstein-Uhlenbeck process is Markov one. It is ergodic (i.e.
possesses unique invariant distribution), if and only if, $\gamma
<0$ and \be\label{logint} \int_{|u|\geq 1}\ln |u|\mu(du)<+\infty;
\ee see \cite{SY84}. Clearly, our standing  assumption
(\ref{esint}) provides (\ref{logint}). Respective {\it stationary
version of the Ornstein-Uhlenbeck process} can be represented as
$$
X_t=\int_{-\infty}^te^{\gamma(t-s)}\,dZ_s, \quad t\in \Re,
$$
which  is clearly of the form (\ref{y}) with $I=\Re$, $f(t,s)=e^{\gamma(t-s)}\1_{s\leq t}$. Conditions on the existence
and smoothness of the distribution  densities for L\'evy driven Ornstein-Uhlenbeck processes were studied in
\cite{Ma04}, \cite{BK08}, \cite{PZ09}, \cite{S10}. In some exceptional stationary cases, the density can be
represented explicitly, see \cite{BNS01}. However, as far as we know, any references concerning general
estimates or a description of the asymptotic behaviour of such a density are not available.

Another example of a process of the type (\ref{y}) is the \emph{fractional L\'evy motion}, defined,
analogously to the fractional Brownian motion, by the stochastic Weyl integral
\be\label{ZHdef}
Z^H(t)={1\over \Gamma(H+1/2)}\int_{\Re}\left[(t-s)_+^{H-1/2}-(-s)^{H-1/2}_+\right]\,dZ_s, \quad t\in \Re,
\ee
where $x_+=\max(x,0)$, and $H\in (0,1)$ is \emph{the Hurst index}; see  \cite{ST94}, \cite{BCI04}, \cite{M06}, \cite{KM10} and  references therein.  In what follows, we will study the asymptotic behaviour of the distribution density of $Z^H(t)$ under the assumption that  $H>1/2$, which is  the so called \emph{long memory case}, see Definition 1.1 in \cite{M06}. Note that in this  case  $Z^H$ is not a Markov process, in contrast to the L\'evy process $Z$, or the L\'evy driven Ornstein-Uhlenbeck process (\ref{sde1}).

Heat kernel estimates for symmetric jump processes were studied
systematically by Barlow, Bass, Chen and Kassman \cite{BBCK09},
Chen, Kim, Kumagai  \cite{CKK08},  \cite{CKK10}, \cite{C09},
Barlow, Grigoryan, Kumagai \cite{BGK09},  Chen, Kumagai
\cite{CK03},  \cite{CK08}, Chen, Kim, Kumagai \cite{CKK09}; see
also Bass and Levin \cite{BL02} for the transition density
estimates for a Markov chain on $\mathbb{Z}^d$. The  approach used
in the papers listed above relies on the paper by  Carlen, Kusuoka
and Stroock \cite{CKS87}. For heat  kernel estimates in domains we
refer to the papers by Bogdan and Jakubowski \cite{BJ07}, Banuelos
and Bogdan \cite{BB04},  Bogdan,  Grzywny, and  Ryznar
\cite{BGR10}, Chen, Kim and Song  \cite{CKSa} -- \cite{CKSd}. Of
course, this list of references is far from  complete.

In particular, heat kernel estimates for symmetric jump processes
on $\mathbb{R}^d$  with jump  kernel $J(x,y)$,  either
bounded both from above and below by
$\frac{1}{|x-y|^{d+\alpha}}1_{|x-y|\leq 1}$, $0<\alpha<2$, $d\geq
1$, or decaying  as $e^{-\gamma |x-y|^\beta}$, $\beta\in [0,\infty)$,
as $|x-y|\to\infty$, are studied in \cite{CKK08} and
\cite{CKK10}, respectively.  Under the particular choice of the
jump kernel $J(x,y)=J(x-y)$, the processes studied in
\cite{CKK08} and \cite{CKK10} become symmetric L\'evy processes. We
postpone to Section \ref{s3} (Example \ref{ex4}) the detailed comparison of the
asymptotic results for the distribution densities of such
processes obtained in \cite{CKK08} and \cite{CKK10}, with the
results obtained by our approach. Here we just mention that
our approach,  based on the the complex analysis technique, can be
applied both for  non-symmetric Markov jump processes, like the
L\'evy driven Ornstein-Uhlenbeck process, and for non-Markov
processes such as  the fractional L\'evy motion.

Let us outline the rest of the paper. Our main result on asymptotic behavior of the distribution densities of  L\'evy
driven stochastic integrals $Y$, including the L\'evy process $Z$ itself, is formulated and proved in Section \ref{s2}.
To simplify the exposition, we give one-sided asymptotics; that is, we formulate the main result for   the distribution
density $p_t(x)$ only for $x\geq 0$. Clearly, one can easily deduce from this result the two-sided asymptotics, assuming
additionally that the L\'evy measure of the process $(-Z)$ satisfies conditions of Theorem \ref{td1}.

Conditions of the main result, Theorem \ref{td1}, are quite abstract, and require an additional analysis in order
to provide a verifiable criteria. For the reader's convenience and to clarify the exposition, we separate such an
analysis in two parts. In Section \ref{s3} we consider an ``individual'' asymptotic behavior of the distribution
density of $Y_t$ with fixed $t$. We formulate an ``individual'' version of Theorem \ref{td1} with verifiable
conditions on the  L\'evy measure $\mu$ and the kernel $f$. These conditions, in particular, reveal  the ``smoothifying''
effect provided by  the kernel $f$: typically, both to provide existence of the distribution density of the L\'evy driven
stochastic integral $Y_t$ and to describe its asymptotic behavior, {fewer} restrictions on the L\'evy measure are required
than in the case of the L\'evy process $Z_t$ itself. {An} illustrative example of such an effect is provided by the fractional L\'evy motion,  where the assumptions on the L\'evy measure are finally reduced to
\be\label{mu0}
\mu(\Re^+)>0.
\ee

 In Section \ref{s4} we establish the asymptotic behavior of the distribution density of $Y_t$, involving both
 state space variable $x$ and   time variable $t$. To shorten the exposition, we restrict ourselves {to} the
 case of a self-similar kernel $f$. The class of the L\'evy driven stochastic integrals with self-similar
 kernels, although not being {the} most general possible, is wide enough to cover the important particular cases
 of the L\'evy process $Z$ itself and the fractional L\'evy motion $Z^H$.  As a corollary of the main result
 of Section \ref{s4} (Theorem \ref{t41}), we obtain asymptotic relation
  \be\label{answ41}
p_t(x)\sim \frac{1}{\sqrt{2\pi t\Kf_Z\left({x\over  t}\right)}} e^{t\Df_Z\left({x\over t}\right)}, \quad
t+x\to \infty,\quad (t,x)\in [t_0, +\infty)\times \ax
\ee
for the distribution density of the L\'evy process $Z$,  and
 \be\label{answ42}
p_t(x)\sim \frac{1}{\sqrt{2\pi t^{2H}\Kf_{Z^H}\left({x\over  t^{H+1/2}}\right)}} e^{t\Df_{Z^H}\left({x\over t^{H+1/2}}\right)}, \quad
t+x\to \infty,\quad (t,x)\in [t_0, +\infty)\times \ax
\ee
for the distribution density of  the fractional L\'evy motion $Z^H$. Here $t_0>0$ is arbitrary, $\Df_Y$, $\Kf_Y$ with $Y=Z$, $Z^H$
are some functions, defined  in terms of the L\'evy measure $\mu$ and the kernel $f$; see Section \ref{s4} below.
Observe that the asymptotic formulae for distribution densities of $Z$ and $Z^H$ possess {the} self-similarity property
in spite of the fact that, in general, the families of these densities are  not self-similar.

Formally, $Z^H$  includes $Z$ as a partial case with $H=1/2$, and
(\ref{answ42}) with $H=1/2$ transforms to (\ref{answ41}). However,
there is a substantial difference between the conditions under
which these asymptotic results are available (see Corollary
\ref{c41}). To get (\ref{answ41}), one should impose some
``regularity'' conditions ($N_1$) and ($C$) together with some
``tail'' conditions ($T_1$) and ($T_2$). To get   (\ref{answ42})
with  $H\in (1/2,1)$,  it is sufficient to claim only ``tail''
conditions and non-degeneracy condition (\ref{mu0}): there is no
need {for} additional ``regularity'' conditions. Such a difference is
caused by the ``smoothifying'' effect provided by  the kernel in
the integral (\ref{ZHdef}).

 Theorem \ref{td1} and Theorem \ref{t41}  describe the asymptotic
behaviour          of the distribution density precisely,
but in an implicit form. In Section \ref{s5} we use these  theorems  in
order to deduce  explicit, although less precise, asymptotic
expressions. In the same section we give another application of Theorem \ref{td1}, and
study the asymptotic
behavior  (as $x\to \infty$ for a fixed $a$) of the ratio
\be\label{c1}
r_a(x)={p(x+a)\over p(x)}
\ee
for the invariant distribution density $p$ of the  Ornstein-Uhlenbeck process. Such a study is
 of particular theoretical interest, since the ratio (\ref{c1})
appears in  the formula for the generator of the dual (i.e.,
time-reversed) process corresponding to  the solution to SDE
(\ref{sde}). Therefore,  knowledge of the asymptotic properties of
(\ref{c1}) would be useful when one is interested in
studying the stationary version of the solution, respective
Dirichlet form etc. For instance, in the forthcoming paper
\cite{Ku09} the estimate given in  Theorem \ref{td4} below is used
substantially in the proof of the spectral gap property for the L\'evy driven
Ornstein-Uhlenbeck process.

Formula (\ref{answ41}) and Theorem \ref{t51} provide a detailed
description of the asymptotic behavior of the distribution
densities of the L\'evy process and the fractional L\'evy motion.
This behavior exhibits two different regimes. In the first regime,
where the ratio $x\over t$ (resp., $x\over t^{H+1/2}$) stays
bounded,  the principal behavior of $p_t(x)$ is determined by the
values of the functions $\Df_Y$, $\Kf_Y$ (with  $Y=Z$ or  $Z^H$)
on a bounded domain. For instance, for any $\varkappa\geq 0$
 \be\label{answ411}
p_t(t\varkappa)\sim \frac{1}{\sqrt{2\pi t\Kf_Z\left(\varkappa\right)}} e^{t\Df_Z\left(\varkappa\right)}, \quad t\to +\infty
\ee
(for the L\'evy process $Z$) and
 \be\label{answ421}
p_t(t^{H+1/2}\varkappa)\sim \frac{1}{\sqrt{2\pi t^{2H}\Kf_{Z^H}\left(\varkappa\right)}} e^{t\Df_{Z^H}\left(\varkappa\right)},  \quad t\to +\infty
\ee
(for the fractional L\'evy motion $Z^H$). In the second regime, where the ratio $x\over t$ (resp., $x\over t^{H+1/2}$) tends to $+\infty$,
 the principal behavior of $p_t(x)$ is determined by the asymptotics of $\Df_Y$, $\Kf_Y$ (with  $Y=Z$ or  $Z^H$) on $+\infty$. Such asymptotics {are} described in Theorem \ref{t41} for two cases: for the L\'evy measure $\mu$ being either ``truncated'' (i.e. supported in a bounded set) or ``exponentially damped'' (i.e. its tail satisfies certain exponential estimate, see (\ref{mu1})).  This description gives some constant $c_*$, determined in terms of the L\'evy measure $\mu$ only (see (\ref{cstar}) and (\ref{cstar1})), such that the statements below hold true (see Corollary \ref{c51} and Corollary \ref{c52}  below).

I. {\itshape Case of the  L\'evy process $Z$.}    For any constants $c_1>c_*$ and $c_2<c_*$ there exists $y=y(c_1,c_2)$ such that for  $x/t>y$, either
\be\label{answ711}
\exp\left(-c_1 x\ln \left(x\over t\right)\right)\leq p_t(x)\leq \exp\left(-c_2 x\ln \left(x\over t\right)\right)
\ee
(if $\mu$ is truncated), or
\be\label{answ712} \exp\left(-c_1x\ln^{\beta-1\over \beta} \left(x\over t\right)\right)\leq p_t(x)\leq \exp\left(-c_2x\ln^{\beta-1\over \beta} \left(x\over t\right)\right)
\ee
(if $\mu$ is is exponentially damped).

II. {\itshape Case of the fractional L\'evy motion $Z^H$.} For any constants $c_1>c_*$ and $c_2<c_*$ there exists $y=y(c_1,c_2)$ such that for  $x/t^{H+1/2}>y$, either
\be\label{answ721}
 \exp\left(-{c_1x \over \Gamma(H+1/2)t^{H-1/2}} \ln \left(x\over t^{H+1/2}\right)\right)\leq p_t(x)\leq \exp\left(-{c_2x\over \Gamma(H+1/2)t^{H-1/2}} \ln \left(x\over t^{H+1/2}\right)\right)
\ee
(if $\mu$ is truncated), or
\be\label{answ722}
\exp\left(-{c_1x \over \Gamma(H+1/2)t^{H-1/2}}\ln^{\beta-1\over \beta} \left(x\over t^{H+1/2}\right)\right)\leq p_t(x)\leq \exp\left(-{c_2x \over \Gamma(H+1/2)t^{H-1/2}} \ln^{\beta-1\over \beta} \left(x\over t^{H+1/2}\right)\right)
\ee
(if $\mu$ is is exponentially damped).

In this paper we restrict ourselves to the case of one-dimensional
processes   in order to make the exposition reasonably short, and
to  give the main results in their most transparent form. These
results have straightforward generalizations to  the
multi-dimensional case; we postpone  the discussion of these
generalizations to a further publication. We also restrict our
considerations of the L\'evy process $Z$  and the fractional
L\'evy motion $Z^H$ to the case where the time variable $t$ is
separated from $0$. The small time estimates require additional
analysis of the local behavior of the L\'evy measure of the noise;
this analysis is performed in the separate article \cite{KK11}.

\section{The main result}\label{s2}
\subsection{Preliminaries}
Everywhere below $Z$ is a L\'evy process and $\psi$ is its
characteristic exponent; that is, (\ref{char}) and (\ref{LKh})
hold.

To  exclude from consideration {the} trivial cases, we assume that $b=0$ and $\mu(\Re)>0$; that is,  $Z$ does not contain a diffusion part, and contains a non-trivial jump part.  Moreover, we assume that $\mu$ satisfies (\ref{mu0}), which is motivated by our intent to analyze the distribution density on the positive half-line. Finally, we assume $Z$
to be centered, which means that the characteristic exponent is
of the form
\begin{equation}\label{LKh1}
 \psi(z)=\int_{\mathbb{R}}
 \left( e^{iuz}-1- iz u \right)\mu(du), \quad z\in \Re.
\end{equation}
This assumption  does not restrict the
generality:  under (\ref{esint}),  the increments of $Z$
have  moments of all orders, therefore the difference between the processes with characteristic exponents
(\ref{LKh}) and (\ref{LKh1}) is given by the explicitly calculable
constant, which clearly does not effect the distributional
properties.

We consider L\'evy driven stochastic integrals of the form
\be\label{int}
Y_t=\int_I f(t,s) dZ_s, \quad t\in \TT,
\ee
where $\TT\subset \Re$ is some set, and $I\subset \Re$ is an interval.  We allow the case
where the interval $I$  belongs not only to the half-line, but to whole  $\mathbb{R}$.
In this case, the process given by   (\ref{int}) is assumed to be well defined on the whole line
$\Re$, and to have independent and stationary distributed increments, with  the
characteristic exponent of the increments still being of the form (\ref{LKh1}).
A standard version of such a process is the so called two-sided L\'evy process
$$
Z_t=\begin{cases}
Z_t^1,&t\geq 0
\\-Z_{-t-}^2, & t<0
\end{cases},
$$
where $Z^1$ and $Z^2$ are two independent copies of a L\'evy process, defined on $\ax$.

We interpret (\ref{int})  as an integral with respect to an infinitely divisible random measure;
for the general theory of such  integrals we refer to \cite{RR89}.
The integral (\ref{int}) is well defined if, and only if,
\be\label{fint}
\int_If^2(t,s)\, ds<+\infty, \quad t\in \TT,
\ee
and in that case its characteristic function admits the representation
\be\label{chint}
 Ee^{izY_t}=\exp\left[\int_I\int_{\Re}\left(e^{izf(t,s)u}-1-izf(t,s)u\right)\,\mu(du)ds\right], \quad z\in\Re, \quad t\in  \TT,
\ee
see Theorem 2.7 from \cite{RR89}. In what follows we assume that  $f$ satisfies (\ref{fint}),  and $f(t,\cdot)$ is bounded for every $t\in \TT$. To exclude the trivial case $Y_t=0$ a.s., we  assume $\int_If^2(t,s)\, ds>0$,  $t\in \TT.$
 We also assume
\be\label{fnon}\int_I(f(t,s)\vee 0)^2\, ds>0, \quad t\in \TT.
\ee
This does not restrict generality since otherwise one can consider $-Y_t$ instead of $Y_t$.

For a Borel set $A\subset \Re$,  denote
$$
\Theta(t,z, A)=\iint_{\{(s,u)\in I\times \Re: \, f(t,s) u\in A\}} (1-\cos(f(t,s) z u))\mu(du)ds, \quad t\in \ax, z\in \Re.
$$
The functions $\Theta(\cdot,\cdot, A)$, with properly chosen sets $A$, will be used below as a tool for studying the properties of distribution densities of  L\'evy driven stochastic integrals $Y$.  One statement of such a type is formulated in the proposition below,  which is in fact  the classic Hartman-Wintner sufficient condition (\cite{HW42}), reformulated in the context of  L\'evy driven stochastic integrals.

As usual, we denote by $C_b^k(\Re)$ the class of function, continuous and bounded together with their derivatives up to order $k$.

\begin{prop}\label{p21} For  given $t\in \TT$, $k\in \mathbb{Z}_+$,  and $|z|$ large enough, let
\be\label{HW1}
\Theta(t,z, \Re)\geq (k+1+\delta)\ln |z|
\ee
with some $\delta>0$.

Then $Y_t$ has a distribution density $p_t$, which belongs to the class $C_b^k(\Re)$.

In particular, if for a given $t\in \TT$
\be\label{HW2}
\Theta(t,z, \Re)\gg \ln |z| \quad \text{as} \quad |z|\to\infty,
\ee
then $Y_t$ has a distribution density $p_t\in C_b^\infty$.
\end{prop}

\proof By (\ref{chint}), condition (\ref{HW1}) implies
\be\label{21}
|Ee^{izY_t}|\leq |z|^{-k-1-\delta} \quad \text{for $|z|$ large enough}.
\ee
Hence the required statement follows by the inversion formula for the Fourier transform.
\endproof

 As usual, we write  $f(\xi)\sim g(\xi)$, $\xi\to\infty$, or
$f(\xi) =o(g(\xi))$, $\xi\to \infty$, if $\lim_{\xi\to\infty}
\frac{f(\xi)}{g(\xi)}=1$ or $\lim_{\xi\to\infty}
\frac{f(\xi)}{g(\xi)}=0$, respectively. We also use the notation
$f(\xi)\ll g(\xi)$, $\xi\to\infty$,  instead of $f(\xi)
=o(g(\xi))$, $\xi\to \infty$, when it is more convenient. The same
conventions are used when functions $f$ and  $g$ depend on $t$
and/or on $x$.

\subsection{The main result: formulation and discussion}
Since $f$ is bounded, the exponential integrability assumption (\ref{esint}) implies  that for $t\in \TT$ the function
$$
\Psi(t,z)=\int_I\int_{\Re}\left(e^{-izf(t,s)u}-1+izf(t,s)u\right)\,\mu(du)ds,  \quad t\in \TT,\quad  z\in \CC,
$$
is well defined and
analytic with respect to  $z$. Denote
$$
H(t,x,z)=ixz+\Psi(t,z),
$$
 and observe  that, assuming (\ref{HW1}), we have
\be\label{ptx3}
 p_t(x)=\frac{1}{2\pi} \int_{\mathbb{R}} e^{H(t,x,z)}\,
dz,\quad x\in \Re,
\ee
which is just the inversion formula for for the characteristic function of $Y_t$, combined with the change of variables $z\mapsto -z$.

Denote
$$
\Mf_{k}(t,\xi)={\prt^k\over \prt \xi^k}\Psi(t,i\xi), \quad
k\geq 1,  \quad \xi\in\mathbb{R}.
$$
Clearly,
$$
\Mf_{k}(t,\xi)=\int_I\int_{\Re}u^kf^k(t,s)e^{\xi f(t,s)u}\,\mu(du)ds={\prt^k\over \prt \xi^k}H(t,x,i\xi), \quad k\geq 2.
$$

Since  $\mu$ and $f(t,\cdot)$ are assumed to be non-degenerate, we have $\Mf_{2}(t,\xi)>0$.
Therefore there exists at most one solution $\xi(t,x)$ to the equation
\be\label{eq}
{\prt\over \prt \xi}H(t,x,i\xi)=0.
\ee
Clearly, for any $t\in \TT$ we have   $\xi(t,0)=0$. Note that
$$
\Mf_{1}(t,\xi)={\prt\over \prt \xi}\Psi(t,x,i\xi)=\int_I\int_{\Re}uf(t,s)\Big(e^{\xi f(t,s)u}-1\Big)\,\mu(du)ds=\int_{I\times \Re}v(e^{\xi v}-1)\mu_{t,f}(dv),
$$
where $\mu_{t,f}$ denotes the image of the measure $\mu(du)ds$ under the mapping
$$
I\times \Re\ni (s,u)\mapsto f(t,s)u\in \Re.
$$
Under the  assumptions (\ref{mu0}) and (\ref{fnon}), which we assume to hold everywhere below, we have  $\mu_{t,f}(\ax)>0$. Therefore  $\Mf_1(t,\xi)\to +\infty$ as $\xi\to+\infty$, which means that  $\xi(t,x)$ is well defined and positive for $x> 0$, and
\be\label{xi}
\xi(t,x)\to +\infty, \quad \text{as} \quad x\to +\infty.
\ee
Note that $z=i\xi(t,x)$  is the unique critical point for  $H(t,x,\cdot)$ on the line $i\Re$.

We put
$$
D(t,x)=H(t,x,i\xi(t,x)),\quad
K(t,x)=\Mf_{2}(t,\xi(t,x))={\prt^2\over \prt
\xi^2}H(t,x,i\xi)\Big|_{\xi=\xi(t,x)}.
$$

In the sequel, we fix $\mathcal{A}\subset \TT\times \Re^+$ and
denote
$$\mathcal{T}=\{t: \exists x\in \Re^+, (t,x)\in \mathcal{A}\},
\quad \mathcal{B}=\{(t,\xi): \exists (t,x)\in \mathcal{A},
(t,\xi)=(t, \xi(t,x))\}.$$ For instance,   if
$\mathcal{A}=\TT'\times \Re^+$ with some  $\TT'\subset \TT$, then
$\mathcal{T}=\TT'$ and $\mathcal{B}=\mathcal{A}$.
% (under the standing  assumptions   (\ref{mu0}) and (\ref{fnon})).

In the following theorem, which represents the main result of the
paper, the function $\theta:\mathcal{T}\to (0,+\infty)$ is assumed to
be bounded away from zero on $\TT$, and the function
$\chi:\mathcal{T}\to (0,+\infty)$ is assumed  to be bounded away
from zero on every set  $\{t:\theta(t)\leq c\}$, $c>0$.
For a particular process $Y$, the choice of the ``scaling''
functions $\theta$ and $\chi$ is determined by the structure of the
kernel $f$, see Section \ref{s4} below.

\begin{theorem}\label{td1} Assume that the following conditions hold true:

($H_1$) $\Mf_{4}(t,\xi)\ll \Mf_{2}^2(t,\xi)$, $\ \theta(t)+\xi\to
\infty$, $(t,\xi)\in\mathcal{B}$.

($H_2$)
$$
\ba
\ln\left(\Big(\chi^{-2}(t)\frac{\Mf_{4}(t,\xi)}{\Mf_{2}(t,\xi)}\Big)
\vee 1\right)&+\ln\left(\Big( \ln
\Mf_{2}(t,\xi)\Big)\vee1\right)\\& \ll \ln
\theta(t)+\chi(t)\xi,\quad \theta(t)+\xi\to \infty,\quad
(t,\xi)\in\mathcal{B}.
\ea
$$

($H_3$) There exist $R>0$ and  $\delta>0$ such that
\be\label{h5}
 \Theta(t,z,\Re^+)\geq (1+\delta)\ln |z|,\quad t\in \mathcal{T},\quad  |z|>R.
\ee

($H_4$) There exists $r>0$ such that for every $\eps>0$,
$$
\inf_{t\in \mathcal{T}, |z|>\eps}\Theta(t,z,[r\chi(t),+\infty))\geq  \theta(t)\Big((\eps\chi(t))^2\wedge 1\Big).
$$

Then for every  $t\in \mathcal{T}$ the law of $Y_t$ has a continuous bounded distribution density $p_t(x)$, and
\begin{equation}\label{answ1}
p_t(x)\sim \frac{1}{\sqrt{2\pi K(t,x)}} e^{D(t,x)}, \quad \theta(t)+x\to \infty, \quad (t,x)\in \mathcal{A}.
\end{equation}

\end{theorem}

\begin{rem}\label{rt1}\textit{(On conditions)}. {The} conditions of Theorem \ref{td1} are rather technical and  abstract. In Sections \ref{s3} and \ref{s4} below we give their more explicit versions, formulated in  terms of the L\'evy measure $\mu$ and the kernel $f$. Note that ($H_1$) and ($H_2$) are, in fact, the assumptions on the  growth of {the}
tails  of the L\'evy measure $\mu$. In addition,  ($H_2$) is balanced with ($H_4$), which  in turn is closely related to the so called {\it  Cramer condition} (see,
for example, \cite{Lu79}, \cite{IL65}, Chapter 3 \S3, and the discussion prior to Lemma \ref{l26} below). Finally, ($H_3$) is a proper uniform version of the Hartman-Wintner condition, see  Proposition \ref{p21}. Clearly, one can consider the stronger version of condition (\ref{h5}) with $k+1+\delta$ instead of $1+\delta$ (cf. (\ref{HW1}) and (\ref{HW2})), and provide the asymptotic relations similar to (\ref{answ1}) for the derivatives of the distribution density $p_t(x)$ up to order $k$.
\end{rem}

\begin{rem}\textit{(On relation (\ref{answ1}))}. 1. Note that the asymptotic relation (\ref{answ1})
 corresponds completely to the standard form of an
asymptotic relation  obtained  by the Laplace method.
Typically, within this method one  can  prove that  the integral
$$
\int_{(a,b)}e^{-F(\lambda, x)}\, dx
$$
is asymptotically equivalent to
\be\label{answl}
\sqrt{\frac{2\pi}{ F''_{xx}(\lambda,
x_\lambda)}}e^{-F(\lambda, x_\lambda)},\quad
x_\lambda:=\mathrm{arg}\min_{x}F(\lambda,x).
\ee
Clearly,  (\ref{answ1}) is  exactly of  the form
(\ref{answl}) with appropriate $F$ and additional normalizer $1/(2\pi)$, which comes from the inverse Fourier transform formula.

2. Our approach is in some sense related to the Large Deviations
Principle (LDP). Namely, if $P_t^l(dx)$ is the probability measure
associated with $Y_t^l:=\frac{1}{l}\sum_{i=1}^l Z_t^i$, where
$\{Z_t^i\}_{i=1}^l$ are independent copies of $(Z_t)_{t\geq 0}$,
then $P_t^l(dx)$ satisfies the LDP  with a {\itshape good rate
function} $\Lambda_t(x):=-D(t,x)$, in the sense that for all measurable subsets $A\subset\mathbb{R}$
\begin{equation*}
    -\inf_{x\in \mathop{\mathrm{interior}}(A)} \Lambda_t(x)
    \leq \liminf_{l\to\infty} \frac{1}{l} \ln P_t^{l}(A)
    \leq \limsup_{l\to\infty} \frac{1}{l} \ln P_t^{l}(A)
    \leq -\inf_{x\in \mathop{\mathrm{closure}}(A)} \Lambda_t(x);
\end{equation*}
 see \cite
{DS}. Moreover, assuming the  exponential integrability condition
\eqref{esint} and  existence of the transition probability
density $p_t(x)$ for $t>t_0$, it is shown in \cite{KS} that
\be\label{KS}
\lim_{l\to\infty} \frac{\ln p_{tl}(lx)}{l}=D(t,x),
\ee
cf. (\ref{answ411}) and (\ref{answ421}).
\end{rem}

\subsection{Proof} Note that $\Theta (t,z,A)$ depends on the set $A$ monotonously. Hence ($H_3$) yields (\ref{HW1}), and therefore (\ref{ptx3}) holds. In what follows, we analyze the right hand side of  (\ref{ptx3}). We divide this   analysis into several steps.

{\it Step 1: changing the integration  contour}. We  prove
that
\begin{equation}
p_t(x)=\frac{1}{2\pi} \int_{i\xi(t,x)+\mathbb{R}} e^{H(t,x,z)}\,
dz=\frac{1}{2\pi} \int_{\mathbb{R}} e^{H(t,x,\eta+i\xi(t,x))}\,
d\eta.  \label{ptx2}
\end{equation}
Recall that we assumed $x\geq 0$, which in turn implies    $\xi(t,x)\geq 0$. Consider the domain
\begin{equation}
G_M:=\Big\{ z\in \mathbb{C}:   \im  z\in [0,\xi(t,x)],\quad  \re
z\in [-M,M], \quad M>0\Big\}.\label{DG}
\end{equation}
The function $H(t,x,z)$ is analytic in $G_M$, hence by the
Cauchy theorem
\begin{equation}
\int_{\partial G_M} e^{H(t,x,z)}\, dz=0.\label{inc}
\end{equation}
Consider the integrals
\begin{equation}
\int_0^1 e^{H(t,x,\pm M + iv\xi(t,x))}dv. \label{c2}
\end{equation}
We have
\begin{equation}
\re H(t,x,\eta+i\xi)= -x\xi - \int_I\int_{\mathbb{R}}
\big(1-e^{f(t,s)\xi u} \cos(f(t,s)\eta u)+f(t,s) \xi u \big)\mu(du)ds
$$
$$
=H(t,x,i\xi)-\int_I\int_{\mathbb{R}} e^{f(t,s)\xi u}(1-\cos(f(t,s) \eta
u))\mu(du)ds, \quad \xi,\eta\in\Re.\label{Re}
\end{equation}
The function $\xi\mapsto H(t,x,i\xi)$ is real-valued, convex, and {attains}
its minimal value at the point $\xi(t,x)$. Then $H(t,x,i\xi)\leq
H(t,x,0)=0$ for $\xi\in[0,\xi(t,x)]$. On the other hand, for every $\xi\geq 0$
$$
\ba
\int_I\int_{\mathbb{R}} e^{f(t,s)\xi u}(1-\cos(f(t,s) \eta
u))\mu(du)ds&\geq \iint_{\{(s,u)\in I\times \Re: \, f(t,s) u>0\} }(1-\cos(f(t,s) \eta
u))\mu(du)ds\\&=\Theta(t,\eta,\ax).
\ea
$$
Therefore
$$
\re H(t,x,\pm M+iv\xi(t,x))\leq -\Theta(t,\pm M,\ax),\quad  v\in[0,1].
$$
Thus,  condition ($H_3$)  implies  that the integrals in \eqref{c2} tend
to 0 as $M\to +\infty$, which together with (\ref{DG}) gives
(\ref{ptx2}).

In what follows we denote
$$R(t,x,\eta)=\re
H(t,x,\eta+i\xi(t,x)),
$$
$$I(t,x,\eta)=\im H(t,x,\eta+i\xi(t,x))=x\eta -t\int_{\mathbb{R}} \big(e^{\xi(t,x) u} \sin(\eta u)-\eta u
\big)\mu(du).
$$
Since a distribution density is real valued, we derive from
\eqref{ptx2}
\begin{equation}
p_t(x)=\frac{1}{2\pi} \int_{\mathbb{R}}
e^{R(t,x,\eta)}\cos(I(t,x,\eta))\, d\eta. \label{aux2}
\end{equation}

Before proceeding further on, let us give a  short description
of the rest of the proof. We will estimate the integral
(\ref{aux2}) using the appropriate version of the
Laplace method (see \cite{Co65} for its description).
In our case {the} application of the Laplace
method meets some difficulties, since  the expression under  the
integral contains two functions $R$ and $I$. Therefore we introduce two intervals
$[-\alpha,\alpha]$ and $[-\beta,\beta]$, on which $R$ and $I$ are
controllable in terms of their Taylor's expansions.  Then we
split the integral into the sum of  integrals over
$\{|\eta|\leq\alpha\}, \{|\eta|\in (\alpha,\beta]\}$, and
$\{|\eta|>\beta\}$, and estimate these integrals separately. As in the standard Laplace method, the
first two integrals are controlled by  using  Taylor expansion
arguments.  For the third
integral, any standard considerations,  like  convexity arguments
from \cite{Fe77}, cannot be applied. Therefore, we use the specific
arguments based on the structure of the functional under
consideration.

{\it Step 2: choosing $\alpha,\beta$.} Following the explanations given above, we split the integral \eqref{aux2} into the  sum
\begin{equation}\label{split}\begin{aligned}
\frac{1}{2\pi} \Big[ \int_{|\eta|\leq \alpha}+\int_{|\eta|\in
(\alpha,\beta]}+\int_{|\eta|>\beta}\Big]
&\Big(e^{R(t,x,\eta+i\xi(t,x))}\cos(I(t,x,\eta+i\xi(t,x)))\,
d\eta\Big)
 \\&=J_1(t,x)+J_2(t,x)+J_3(t,x),\end{aligned}
\end{equation}
 where  $\alpha\equiv \alpha(t,x)$
and $\beta\equiv \beta(t,x)$ are auxiliary functions. The
function $\beta$ is defined by \be\label{beta}
\beta(t,x)=\sqrt{\frac{\Mf_{2}(t,\xi(t,x))}{\Mf_{4}(t,\xi(t,x))}}.
\ee Our aim in this step is to construct the function $\alpha$
in such a way that
\begin{equation}0<\alpha(t,x)\leq \beta(t,x), \quad (t,x)\in \Af,\label{a01}
\end{equation}
\begin{equation}\ba
\frac{1}{\Mf_{2}(t,\xi(t,x))}\ll \alpha^2(t,x)&\ll
\frac{\Mf_{2}(t,\xi(t,x))}{\Mf_{4}(t,\xi(t,x))}, \\
\alpha^3(t,x)&\ll \frac{1}{\Mf_{3}(t,\xi(t,x))}, \quad
\theta(t)+x\to \infty, \quad (t,x)\in \Af.\ea \label{a1}
\end{equation}
By the Cauchy inequality and condition ($H_1$), we
have
$$
\Mf_{3}^2(t,\xi)\leq \Mf_{2}(t,\xi)\Mf_{4}(t,\xi)\ll
\Mf_{2}^3(t,\xi),\quad \theta(t)+\xi\to \infty, \quad (t,\xi)\in
\mathcal{B}.
$$
Hence, there exists a function $\kappa=\kappa(t,\xi)$, such that
\begin{equation}\ba
1\ll \kappa(t,\xi), \quad & \kappa(t,\xi)\ll \Mf_{2}(t,\xi)
\Mf_{4}^{-1/2}(t,\xi), \\ & \kappa(t,\xi)\ll
\Mf_{2}^{1/2}(t,\xi)\Mf_{3}^{-1/3}(t,\xi),\quad
\theta(t)+\xi\to \infty, \quad (t,\xi)\in \mathcal{B}.\ea
\label{mu}
\end{equation}
Without loss of generality, we can assume the function $\kappa$ to be locally bounded. Then we put
\begin{equation}
\alpha(t,x)=c \kappa(t,\xi(t,x)) \Mf_{2}^{-1/2}(t,\xi(t,x))\label{ab}
\end{equation}
with some constant $c>0$. By  (\ref{mu}) and  (\ref{xi}), we have (\ref{a1}). Since $\kappa$ is locally bounded, the constant $c$ can be chosen small enough to provide (\ref{a01}).

 {\it Step 3: estimating $J_1(t,x)$ in
(\ref{split})}. {A} straightforward computation shows that
\begin{equation}\label{estR}
\frac{\partial }{\partial \eta} R(t,x,
\eta)\big|_{\eta=0}=\frac{\partial^3  }{\partial \eta^3} R(t,x,
\eta)\big|_{\eta=0}=0, \quad \frac{\partial^2 }{\partial \eta^2}
R(t,x, \eta)\big|_{\eta=0}=-\Mf_{2}(t,\xi(t,x)),
$$
$$ \Big|\frac{\partial^4 }{\partial \eta^4}
R(t,x, \eta)\Big|=\left|\int_I\int_{\Re}u^4f^4(t,s)e^{\xi f(t,s)u}\cos\Big(\eta f(t,s)u\Big)\,\mu(du)ds\right| \leq \Mf_{4}(t,\xi(t,x)),\quad \eta\in \Re,
\end{equation}
which gives
\begin{equation}
-\Mf_{2}(t,\xi(t,x))-\frac{\eta^2}{2} \Mf_{4}(t,\xi(t,x))\leq
\frac{\partial^2}{\partial \eta^2} R(t,x, \eta)\leq -
\Mf_{2}(t,\xi(t,x))+\frac{\eta^2}{2} \Mf_{4}(t,\xi(t,x))  \label{r2est}
\end{equation}
for all $\eta\in\mathbb{R}$. Therefore  by the estimate for $\alpha^2$ in \eqref{a1} we get
\begin{equation}\label{Req}
\ba
&\sup_{|\eta|\leq \alpha} \frac{\partial^2}{\partial \eta^2}
R(t,x, \eta)\sim -\Mf_{2}(t,\xi(t,x)), \\ &\inf_{|\eta|\leq
\alpha} \frac{\partial^2}{\partial \eta^2} R(t,x, \eta)\sim
-\Mf_{2}(t,\xi(t,x)), \quad \theta(t)+x\to \infty, \quad
(t,x)\in \Af.
\ea
\end{equation}
Next, similarly to \eqref{estR} we get
$$
I(t,x,\eta)\big|_{\eta=0} =\frac{\partial}{\partial
\eta}I(t,x,\eta)\big|_{\eta=0} =\frac{\partial^2}{\partial
\eta^2 }I(t,x,\eta)\big|_{\eta=0}=0, \quad \big|\frac{\partial^3
}{\partial \eta^3} I(t,x, \eta)\big|\leq \Mf_{3}(t,\xi(t,x)).
$$
Note that the  equality for ${\prt\over \prt\eta}I$  holds true
because $z=i\xi(t,x)$ is a critical point for $H(t,x,\cdot)$.
Hence the estimate for $\alpha^3$ in \eqref{a1} implies
\begin{equation}\label{Ieq}
\sup_{|\eta|\leq \alpha} |I(t,x,\eta)|\to 0, \quad
\theta(t)+x\to \infty, \quad (t,x)\in \Af.
\end{equation}

Recall that $K(t,x)\equiv\Mf_{2}(t,\xi(t,x))$ and
$D(t,x)\equiv H(t,x, i\xi(t,x))=R(t,x,0)$. From (\ref{Req})
and (\ref{Ieq}) we get
\begin{equation}
\begin{aligned}
\int_{|\eta|\leq \alpha}e^{R(t,x,\eta)}\cos I(t,x,\eta)d\eta &
\sim e^{R(t,x,0)} \int_{|\eta|\leq \alpha}e^{-\frac{K(t,x)
\eta^2}{2}}d\eta
\\&
= \sqrt{\frac{2\pi}{K(t,x)}} e^{R(t,x,0)}
\int_{|\eta|\leq \sqrt{K(t,x)}\alpha}
\frac{e^{-\frac{|\eta|^2}{2}}}{\sqrt{2\pi}}d\eta\\
&  \sim \sqrt{\frac{2\pi}{K(t,x)}} e^{D(t,x)}, \quad
\theta(t)+x\to \infty, \quad (t,x)\in \Af,
\end{aligned}
\end{equation}
where in the last relation we  used the lower estimate for
$\alpha$ in \eqref{a1}. Thus,
\begin{equation}
J_1(t,x)\sim \frac{1}{\sqrt{2\pi K(t,x)}} e^{D(t,x)}, \quad
\theta(t)+x\to \infty, \quad (t,x)\in \Af.\label{i1}
\end{equation}

{\it Step 4: proving that $J_2(t,x)$ in (\ref{split}) is
negligible.} On the set $\{|\eta|\leq \beta\}$, the function $R$
is controlled by its Taylor expansion. Hence  for
the integral $J_2(t,x)$  we can apply standard arguments of the Laplace
method.

 By \eqref{r2est} we have for $|\eta|\leq \beta$
$$
R(t,x,\eta)\leq R(t,x,0)-\frac{1}{4} \Mf_{2}(t,\xi(t,x))
\eta^2,
$$
which, together with  the lower estimate for $\alpha$ in
\eqref{a1}, gives
\begin{equation}
\begin{split}
|J_2(t,x)|& \leq \int_{|y|\in (\alpha,\beta]}
e^{R(t,x,\eta)}d\eta \leq e^{R(t,x,0)} \int_{|y|>\alpha}
e^{-\frac{\Mf_{2}(t,\xi(t,x)) \eta^2}{4}}d\eta
\\&=\frac{e^{D(t,x)}}{\sqrt{K(t,x)}}
\int_{|y|>\alpha\sqrt{K(t,x)}} e^{-\frac{\eta^2}{2}}\, d\eta \ll
J_1(t,x), \quad \theta(t)+x\to \infty, \quad (t,x)\in \Af.
\end{split}
\end{equation}

{\it Step 5: proving that $J_3(t,x)$ in (\ref{split}) is
negligible.}  By  \eqref{Re},
    \begin{align*}
        |J_3(t,x)|&\leq \frac{1}{2\pi} \int_{|\eta|>\beta} e^{R(t,x,\eta)}d\eta
        \\&
        \leq \frac{1}{2\pi} e^{D(t,x)} \int_{|\eta|>\beta} \exp\left\{-\int_I\int_{\mathbb{R}} e^{f(t,s)\xi u}(1-\cos(f(t,s) \eta
u))\mu(du)ds\right\}\, d\eta.
    \end{align*}
Therefore, by  \eqref{i1},  to prove $J_3(t,x)\ll J_1(t,x)$ we need to check that
\begin{equation}
\int_{|\eta|>\beta} e^{-\Delta(t,x,\eta)}d\eta\ll K^{-1/2}(t,x)=
\Mf_{2}^{-1/2}(t,\xi(t,x)),\quad \theta(t)+x\to \infty, \quad
(t,x)\in \Af,\label{bb}
\end{equation}
where
    $$
     \Delta(t,x,\eta)=\int_I\int_{\mathbb{R}} e^{f(t,s)\xi(t,x) u}(1-\cos(f(t,s) \eta
u))\mu(du)ds.
    $$
 Recall that  $\xi(t,x)\geq 0$. Then for any $\sigma\in (0,1)$ we have
\begin{equation}\begin{aligned}
\Delta(t,x,\eta)&\geq\iint_{\{(s,u)\in I\times \Re: f(t,s) u\geq 0\}}  e^{f(t,s)\xi(t,x) u}(1-\cos(f(t,s) \eta
u))\mu(du)ds\\&\geq (1-\sigma) \Theta(t,\eta,\ax)+\sigma e^{r\chi(t)\xi(t,x)} \Theta(t,\eta,[r\chi(t),+\infty)).\end{aligned}
\end{equation}

 Condition ($H_3$), combined with the trivial observation that $\Theta(t,\eta,\ax)$ is non-negative, yields
\begin{equation}
\sup_{t\in \mathcal{T}}\int_{\mathbb{R}} e^{- (1-\sigma)\Theta(t,\eta,\ax)}\,
d\eta<\infty,\label{d1}
\end{equation}
provided  that $\sigma$ is chosen such  that $(1-\sigma)(1+\delta)>1$. Applying condition ($H_4$) with  $\eps=\beta(t,x)$ gives
 $$
 e^{r\chi(t)\xi(t,x)} \Theta(t,\eta,[r\chi(t),+\infty))\geq c e^{r\chi(t)\xi(t,x)}  \theta(t)\Big((\beta(t,x)\chi(t))^2\wedge 1\Big), \quad |\eta|\geq \beta(t,x).
 $$
Thus, in view of
\eqref{d1},  to show \eqref{bb} it is enough  to prove
\begin{equation}
\exp\left[-\sigma e^{r\chi(t)\xi(t,x)}
\theta(t)\Big((\beta(t,x)\chi(t))^2\wedge 1\Big)\right] \ll
\Mf_{2}^{-1/2}(t,\xi(t,x)), \quad \theta(t)+x\to \infty, \quad
(t,x)\in \Af \label{bbb}
\end{equation}
for every $\sigma>0$.  By the definition (\ref{beta}) of
$\beta(t,x)$, we have
\be\label{beta1}
\Big((\beta(t,x)\chi(t))^2\wedge
1\Big)=\left(\Big(\chi^{-2}(t)\frac{\Mf_{4}(t,\xi)}{\Mf_{2}(t,\xi)}\Big)
\vee 1\right)^{-1}.
\ee
Condition $(H_2)$ and  the
assumptions on $\theta$ and $\chi$, imposed prior  to Theorem
\ref{td1}, yield for any $\sigma>0$
$$
\left(\Big(\chi^{-2}(t)\frac{\Mf_{4}(t,\xi)}{\Mf_{2}(t,\xi)}\Big) \vee 1\right)\ln \Mf_{2}(t,\xi)\ll \sigma\theta(t)e^{r\chi(t)\xi},\quad
\theta(t)+\xi\to \infty, \quad (t,\xi)\in \mathcal{B}.
$$
This relation, combined with (\ref{beta1}), (\ref{xi}), and the relation $\xi(t,x)\geq 0$, yields
$$
\ln \Mf_{2}(t,\xi(t,x))\ll \sigma \theta(t) e^{r\chi(t)\xi(t,x)}
\Big((\beta(t,x)\chi(t))^2\wedge 1\Big),\quad \theta(t)+x\to
\infty, \quad (t,x)\in \Af,
$$
which in turn implies (\ref{bbb}) and
completes the proof  of (\ref{bb}).

We have proved
 $$
 J_2(t,x)\ll J_1(t,x),\quad J_3(t,x)\ll
J_1(t,x).
$$   By
\eqref{i1} we get the statement of the theorem.  \hfill $\Box$

\section{Explicit conditions: fixed time setting}\label{s3}
 Our further aim is to give explicit and tractable  sufficient conditions which provide assumptions ($H_1$) -- ($H_4$) of Theorem \ref{td1}.  In this section  we consider  the case where the time variable is fixed. Therefore everywhere  below in this section we assume
$$
\TT=\{t\}, \quad \mathcal{A}=\mathcal{B}=\{t\}\times \Re, \quad \mathcal{T}=\{t\}.
$$
We skip the variable $t$ in the notation and write, for instance,  $Y$, $f(s)$, $\Mf_k(\xi)$ instead of $Y_t,$ $f(t,s),$ $\Mf_k(t,\xi)$, respectively.

In the fixed time setting the assumptions ($H_1$) -- ($H_4$) look more simple: in particular,  functions $\theta(t)$ and $\chi(t)$ degenerate to some constants $\theta$ and $\chi$. Therefore it is appropriate to introduce the set of conditions which will be useful later on.

($\hat H_1$) $\Mf_{4}(\xi)\ll \Mf_{2}^2(\xi)$, $\xi\to +\infty.$

($\hat H_2$)  $\ln\left(\Big(\frac{\Mf_{4}(\xi)}{\Mf_{2}(\xi)}\Big) \vee 1\right)+\ln \ln
\Mf_{2}(\xi)\ll  \xi$,  $\xi\to +\infty$.

($\hat H_3$) There exist $R>0$, $\delta>0$ such that
\be\label{h51}
\Theta(z,\Re^+)\geq (1+\delta)\ln |z|,\quad |z|>R.
\ee

($\hat H_4$) There exist $q>0$ and $\vartheta>0$ such that for every $\eps>0$
$$
 \inf_{|z|>\eps}\Theta(z,[q,+\infty))\geq \vartheta \Big(\eps^2\wedge 1\Big).
$$

One can easily see that in the fixed time setting conditions ($H_1$) -- ($H_4$)  are equivalent to ($\hat H_1$) -- ($\hat H_4$).
Indeed,   the constants $\theta>0$ and $\chi>0$, which come, respectively,  from the functions $\theta(t), \chi(t)$, are suppressed in
($H_2$) by the term $\xi$. In ($H_4$), the constant $\chi$ can be eliminated by a proper change of the constants $r$ and $\theta$;
we denote these new  constants by $q$ and $\vartheta$.
%instead of $r, \theta$.

Clearly, $Y$ is  infinitely divisible  with the L\'evy measure
$$
\mu_f(A)=\iint_{I\times \Re}\1_{uf(s)\in A}\mu(du)ds.
$$
In what follows we demonstrate that  conditions ($\hat H_1$) -- ($\hat H_4$), which are  in fact the assumptions on $\mu_f$, can be verified efficiently in the terms of the kernel $f$ and the initial L\'evy measure $\mu$.

\subsection{Assumptions ($\hat H_1$) and ($\hat H_2$)}
 Observe that ($\hat H_1$) and ($\hat H_2$) control the growth rate of the ``tails'' of $\mu_f$.  The following two {lemmas} show that  these assumptions can be
verified in the terms of similar ``tail'' conditions imposed on $\mu$. Denote
\be\label{m1}
M_1(\xi)=\int_{\Re}u(e^{\xi u}-1) \mu(du), \quad
M_k(\xi)=\int_{\Re}u^ke^{\xi u} \mu(du), \quad k\geq 2.
\ee
Clearly,
\be\label{m2}
\Mf_{k}(\xi)=\int_I f^k(s) M_k(f(s) \xi)\, ds, \quad k\geq 1.
\ee

\begin{lem}\label{l22} Assume

 (${T}_1$) there exists $\gamma\in(0,1)$ such that
$M_4(\xi)\ll M_2^2(\gamma \xi)$, $\xi\to+\infty$.

Then ($\hat H_1$) holds true.
\end{lem}
\proof  Under the assumption (\ref{mu0}) we have $M_k(\xi)\to
+\infty$, $\xi\to +\infty$, $k\geq 1$. In addition, by H\"older
inequality,   for every $a\in (0,1)$ and  $k\geq 2$,
$$
M_k(a\xi)\leq [M_k(\xi)]^a\left[\int_{\Re}u^k\mu(du)\right]^{1-a}
$$
implying
\be\label{22}
M_k(a\xi)\ll M_k(\xi), \quad \xi \to+\infty.
\ee
Denote $F=\mathrm{esssup}_{s\in I} f(s)$ and $I_{f,\gamma}=\{s: f(s)\geq \gamma F \}$. Recall that $f$ is assumed to be bounded, which together with (\ref{fnon}) yields $F\in (0,+\infty)$.  Since  $f^2$ is integrable on $I$ and $f$ is bounded, $f^k, k\geq 3$ is integrable as well. Then (\ref{22}) yields
\be\label{24}
\Mf_k(\xi)=\left(\int_{I_{f,\gamma}}f^k(s)M_k(f(s)\xi)\,ds\right)\Big[1+o(1)\Big], \quad \xi \to +\infty.
\ee
Since $M_2$ is convex and $M_2(\xi)\to +\infty$, as $\xi\to +\infty$, there exists $\xi_0$ such that $M_2$ is increasing on $[\xi_0,+\infty)$. Then for $\xi>\xi_0\gamma^{-1}F^{-1}$ we have
$$\ba
\left(\int_{I_{f,\gamma}}f^2(s)M_2(f(s)\xi)\,ds\right)^2&=\int_{I_{f,\gamma}}\int_{I_{f,\gamma}}f^2(s_1)f^2(s_2)M_2(f(s_1)\xi)M_2(f(s_2)\xi)\,ds_1ds_2
\\& \geq M_2^2(\gamma F\xi) \left(\int_{I_{f,\gamma}}f^2(s)\,ds\right)^2.\ea
$$
Similarly, for sufficiently large $\xi$ we have
$$
\int_{I_{f,\gamma}}f^4(s)M_4(f(s)\xi)\,ds\leq M_4(F\xi)\int_{I_{f,\gamma}}f^4(s)\,ds.
$$
These relations,  together with (\ref{24}) and ($T_1$), imply  ($\hat H_1$).
\endproof

\begin{lem}\label{l23} Assume

 (${T}_2$) $\ln\Big(\frac{M_4(\xi)}{M_2(\xi)}\vee 1\Big) +\ln \ln
M_2(\xi)\ll \xi$, $\xi\to+\infty$.

Then ($\hat H_2$) holds true.
\end{lem}

\proof Fix an arbitrary  $\gamma \in(0,1)$.  For $\xi$ large enough, we have  by (\ref{24})
$$
\Mf_2(\xi)\sim \int_{I_{f,\gamma}}f^2(s)M_2(f(s)\xi)\,ds\leq M_2(F \xi) \left(\int_{I_{f,\gamma}}f^2(s)\,ds\right), \quad \xi \to+\infty,
$$
which together with ($T_2$) gives
\be\label{25}
\ln\ln \Mf_2(\xi)\ll \xi, \quad \xi \to +\infty.
\ee
On the other hand, ($T_2$) implies  that  for every $\eps>0$ and $\xi$ large enough,
$$
M_4(\xi)\leq e^{\eps \xi}M_2(\xi).
$$
Then for $\xi$ large enough  we have by (\ref{24})
$$
\ba \Mf_4(\xi)&\sim \int_{I_{f,\gamma}}f^4(s)M_4(f(s)\xi)\,ds\leq
e^{\eps F\xi} \int_{I_{f,\gamma}}f^4(s)M_2(f(s)\xi)\,ds\\&\leq
F^2 e^{\eps F\xi} \int_{I_{f,\gamma}}f^2(s)M_2(f(s)\xi)\,ds\sim
F^2 e^{\eps F\xi}  \Mf_2(\xi).  \ea
$$
Consequently,
$$
\limsup_{\xi \to +\infty} \xi^{-1} \ln\left({\Mf_4(\xi)\over \Mf_2(\xi)}\right)\leq \eps F.
$$
Since  $\eps>0$ is arbitrary, this relation combined with (\ref{25}) implies ($\hat H_2$).
\endproof

\subsection{Assumptions ($\hat H_3$) and ($\hat H_4$)} It will be convenient to consider, together with the assumption ($\hat H_3$), its  stronger  version
$$
\Theta(z,\Re^+)\gg \ln |z|,\quad z\to \infty.\leqno(\hat H_3^s)
$$
To proceed with the assumption ($\hat H_3^s$),  we introduce several conditions  on the kernel $f$.

($F_1$) $\int_I (f(s)\vee 0)^2\, ds>0$.

($F_2$) On some  interval $[a,b]\subset I$,  the function $f$ is positive and has a continuous non-zero derivative.

($F_3$) On some interval $(-\infty, b]\subset I$, the function $f$ is positive, convex, and has at most exponential decay at $-\infty$;  that is, there exists $\gamma>0$ such that
\be\label{exp}
\lim_{s\to -\infty}e^{-\gamma s}f(s)=+\infty.
\ee

($F_4$) On some interval $(-\infty, b)\subset I$, the function $f$ is positive, convex, and has a subexponential decay at $-\infty$; that is, (\ref{exp}) holds true for every $\gamma>0$.

Note that when $i$ increases from $1$ to $4$, the respective conditions $(F_i)$ become  stronger. Condition  $(F_1)$ is just our standing non-degeneracy assumption (\ref{fnon}), listed here for further reference convenience.  Conditions $(F_1)$ -- $(F_4)$ are well designed to handle the particularly interesting classes of processes, mentioned in the Introduction. Namely,
\begin{enumerate}
\item for the L\'evy process $Z$, one has $f(t,s)=\1_{[0,t]}(s)$, which satisfies $(F_1)$ for every $t>0$;
\item for the non-stationary version of a L\'evy driven Ornstein-Uhlenbeck process, one has $f(t,s)=e^{\gamma(t-s)}\1_{[0,t]}(s)$, which satisfies $(F_2)$ for every $t>0$;
\item for the stationary version of a L\'evy driven Ornstein-Uhlenbeck process, one has $f(s)=e^{\gamma s}\1_{(-\infty, 0]}(s)$, which satisfies $(F_3)$ with $b=0$;
\item for the fractional  L\'evy motion, one has $f(t,s)={1\over \Gamma(H+1/2)}\left[(t-s)_+^{H-1/2}-(-s)^{H-1/2}_+\right]$, which for every $t>0$ satisfies $(F_4)$  with $b=0$.
\end{enumerate}

Recall several conditions which appeared in the literature in the context of the problem of studying local properties of infinitely divisible distributions.

A L\'evy measure $\nu$ on $\Re$ is said to satisfy  the {\it
Hartman-Wintner condition} (\cite{HW42}), if
\begin{equation}\label{HW}
\int_{\Re}(1-\cos zu)\nu(du) \gg \ln |z|, \quad z\to \infty.
\end{equation}
Clearly, ($\hat H_3^s$) is exactly the assumption on  $\mu_f$, restricted to $\Re^+$, to satisfy the Hartman-Wintner condition.

An elementary inequality
\be\label{cos}
c x^2\1_{|x|\leq 1}\leq 1-\cos x\leq x^2\wedge 1, \quad x\in \Re
\ee
(where $c>0$ is some constant) provides the following pair of conditions, sufficient and necessary for (\ref{HW}), respectively:
\begin{equation}\label{Kal}
\int_{|u|\leq |z|^{-1}} (uz)^2 \nu(du) \gg  \ln |z|, \quad
|z|\to \infty; \end{equation}
\begin{equation}\label{Kul}
\int_{\Re} [(uz)^2\wedge 1] \nu(du) \gg \ln |z|, \quad |z|\to
\infty.
\end{equation}
Condition (\ref{Kal}) was introduced in \cite{Ka81}, and is  called  the  {\it Kallenberg condition}. Condition (\ref{Kul}) was introduced in \cite{Ku06}, where it was proved to be necessary for the
 existence of a {\it bounded} transition probability density of the solution  to a (not necessarily linear) L\'evy driven
 SDE.  At  the same time, for an Ornstein-Uhlenbeck process (\ref{sde}) with
 non-trivial drift ($\gamma\not=0$), this condition is sufficient
 for  the existence of $C^\infty$  distribution density (\cite{BK08}).
 Thus, for the non-stationary version of the Ornstein-Uhlenbeck process (\ref{sde}), condition
 (\ref{Kul}) is a criterion.

We denote by  $\mu_+$ the restriction of $\mu$ to $\Re^+$, and formulate the following  set of ``non-degeneracy'' conditions on the measure $\mu$.

($N_1$)  $\mu_+$  satisfies  (\ref{Kal}).

($N_2$) $\mu_+$  satisfies (\ref{Kul}).

($N_3$) $\mu(\Re^+)=+\infty$.

($N_4$) $\mu(\Re^+)>0$.

Note that when $i$ increases  from $i=1$ to $i=4$, the respective conditions $(N_i)$ become more mild; $(N_4)$ is just our fixed non-degeneracy assumption (\ref{mu0}), listed here for further reference convenience.

\begin{lem}\label{l24}  Assume for some  $i=1,\dots, 4$ conditions $(N_i)$ and $(F_i)$ hold.

Then ($\hat H_3^s$) holds true.
\end{lem}
\proof {\it Case $i=1$.}  From the positivity of $1-\cos x$ and the first inequality in (\ref{cos}), it follows that
$$\ba
\Theta(z,\ax)&=\iint_{(s,u):uf(s)>0}(1-\cos(uf(s)z))\mu(du)ds \\& \geq\iint_{(s,u):u>0, 0<uf(s)<1/|z|}(1-\cos(uf(s)z))\mu(du)ds\\& \geq c\iint_{(s,u):u>0, 0<uf(s)<1/|z|}u^2f^2(s)z^2\mu(du)ds\geq cz^2\left[\int_{I}f_+(s)^2\,ds\right]\int_{(0,(F|z|)^{-1})}u^2\mu(du),
\ea
$$
here we keep the notation $F=\mathrm{esssup}_{s\in I} f(s)$. Combined with (\ref{Kal}) for $\mu_+$, the estimates above  provide
 ($\hat H_3^s$).

{\it Case $i=2$.} Since $f$ is positive on $[a,b]$, we have
$$
\Theta(z,\ax)\geq \iint_{(a,b)\times \ax}(1-\cos(uf(s)z))\,\mu(du)ds.
$$
Let us show  that
\be\label{ineq}
\int_{a}^b(1-\cos(xf(s)))\, ds \geq c(x^2\wedge 1)
\ee
holds true with some  constant $c>0$, which would imply  ($\hat H_3^s$) provided that  the assumption (\ref{Kul}) is satisfied. Consider the function
$$
\Upsilon(x)=\int_{a}^b(1-\cos(xf(s)))\, ds.
$$
Clearly, $\Upsilon(x)\sim c_1x^2$ as $x\to 0$, with $c_1=(1/2)\int_{a}^bf^2(s)\, ds>0$. Further, one can write
\be\label{ups}
\Upsilon(x)=\int_{f(a)}^{f(b)}(1-\cos(xv))g'(v)\, dv,
\ee
 where $g:= f^{-1}$.   By our assumptions on  $f$ we have   $g\in C^1$, which implies $\Upsilon(x)>0$  for every $x\not =0$.  Finally, by the Riemann-Lebesgue lemma,
$$
\int_{f(a)}^{f(b)}\cos(xv)g'(v)\, dv\to 0, \quad x\to \infty,
$$
which implies $\lim_{x\to \infty}\Upsilon(x)>0$ and completes the proof of (\ref{ineq}).

{\it Cases $i=3$ and $i=4$.}  We show that the inequality
\be\label{ineq1}
\int_{-\infty}^b (1-\cos(xf(s)))\, ds \geq c \ln |x|
\ee
holds true (i) \emph{for some} $c>0$ and $|x|$ large enough provided that  $f$ satisfies ($F_3$); (ii) \emph{for every} $c>0$ and $|x|$ large enough provided that $f$ satisfies ($F_4$). Keeping the notation $g$ for the inverse function for $f$, we have
$$
\Upsilon(x):=\int_{-\infty}^b(1-\cos(xf(s)))\, ds=\int_0^{f(b)}(1-\cos (xv))g'(v)\, dv.
$$
Since $f$ is convex, $f'$ is non-decreasing. In addition, $f$  itself is increasing: this follows from the convexity, condition (\ref{exp}), and the fact that $f(s)\to 0$ as $ s\to -\infty$ (which comes from the integrability of $f^2(s)$). Therefore, $g'(v)=[f'(g(v))]^{-1}$ is positive and non-increasing.

By positivity of $g'$,
\be\label{ineq3}
\Upsilon(x)\geq \int_{\pi/(2|x|)}^{f(b)}(1-\cos (xv))g'(v)\, dv
\ee
when $\pi/(2|x|)\leq f(b)$. Denote $I_k:=\left[{(2k-1)\pi\over 2|x|},{(2k+1)\pi\over 2|x|}\right], k\geq 1$.
Then, since  $g'$ is positive,
$$
(-1)^k\int_{I_k}\cos (xv)g'(v)\, dv > 0
$$
for every $k\geq 1$, and since $g'$ is non-increasing, we have
\be\label{ineq2}
\int_{I_{k-1}}\cos (xv)g'(v)\, dv+\int_{I_k}\cos (xv)g'(v)\, dv\leq 0
\ee
for every even $k\geq 2$.  Note that, on the axis $[0,+\infty)$, the ``negative'' interval $I_{k-1}$ is located to the left from the ``positive'' interval $I_k$. Then, for any $A>0$, inequality  (\ref{ineq2}) still holds true with $I_{k-1}$ and $I_k$ replaced, respectively, by $I_{k-1}\cap [0,A]$ and $I_{k}\cap [0,A]$. Consequently,   for any $A\geq \pi/(2|x|)$
$$\ba
\int_{\pi/(2|x|)}^{A}\cos (xv)g'(v)\, dv&=\sum_{k=1}^\infty \int_{I_k\cap [0,A]}\cos (xv)g'(v)\, dv\\&=\sum_{m=1}^\infty\left(
\int_{I_{2m-1}\cap [0,A]}\cos (xv)g'(v)\, dv+\int_{I_{2m}\cap [0,A]}\cos (xv)g'(v)\, dv\right)\leq 0.\ea
$$
Therefore we obtain  by (\ref{ineq3})
\be\label{ineq4}
\Upsilon(x)\geq \int_{\pi/(2|x|)}^{f(b)}g'(v)\,dv=g(f(b))-g(\pi/(2|x|))=b-g(\pi/(2|x|))
\ee
for $|x|$ large enough. It follows from (\ref{exp}) that
$$
\rho:=\liminf_{v\to 0}\left(-{g(v)\over \ln (1/v) }\right)
$$
is positive when $f$ satisfies ($F_3$), and equals to $+\infty$ when $f$ satisfies ($F_4$). Combined with (\ref{ineq4}), this yields (\ref{ineq1}).

Now we can complete the proof. In the case $i=3$, take $c>0$ and $Q>0$ such that (\ref{ineq1}) holds true for $|x|\geq Q$. Since  $\mu(\ax)=+\infty$, there exists $q>0$ such that $\mu([q,+\infty))\geq (1+\delta)c^{-1}$. Then (\ref{ineq1}) with $x=uz$ implies
\be\label{Th}
\Theta(z,\ax)\geq \int_{[q,+\infty)}\left(\int_{-\infty}^b (1-\cos(uf(s)z))\,ds\right)\mu(ds)\geq c\mu([q,+\infty))\ln(|qz|), \quad |z|\geq q^{-1}Q,
\ee
 which provides ($\hat H_3^s$) because  $\ln(|qz|)\sim \ln|z|$, $|z|\to \infty$.

In the case $i=4$, {the} assumption $\mu(\ax)>0 $ implies the existence of  $q>0$ for which  $\mu([q,+\infty))>0$. Take $c$ satisfying $c\mu([q,+\infty))>(1+\delta)$, and let $Q>0$ be such that (\ref{ineq1}) holds true with this $c$ and $|x|\geq Q$. Then (\ref{Th}) holds true as well, which provides  ($\hat H_3^s$).
\endproof

Lemma \ref{l24} shows that the kernel $f$ is ``smoothifying'' in the  following sense: when $f$ satisfies some additional assumption like $(F_2)$ -- ($F_4$), the Hartman-Wintner type condition ($\hat H_3^s$) holds true under milder assumptions on the L\'evy measure of the noise.  The following lemma shows that such ``smoothifying'' effect concerns the condition  ($\hat H_4$), as well.

\begin{lem}\label{l25} Under the assumption   (\ref{mu0}) assume additionally that the function $f$ satisfies ($F_2$).

Then  ($\hat H_4$) holds true for $q>0$ small enough.

\end{lem}
\proof
Similarly to the proof of Lemma \ref{l24}, case $i=2$, we assume that $f$ is positive on $[a,b]$. Take $\rho>0$  such that  $\mu([\rho,+\infty))>0$. Then, for  $0<q<\rho\min_{s\in(a,b)}f(s)$, we have by (\ref{ineq})
$$\ba
\Theta(z,[q,+\infty))&\geq \int_{u\geq \rho}\int_{a}^b (1-\cos(uf(s)z))ds\mu(du)\\& \geq c \int_{u\geq \rho}\Big((uz)^2\wedge 1\Big) \mu(du)\geq c\mu([\rho,+\infty))\Big((\rho z)^2\wedge 1\Big),\ea
$$
which implies  the required estimate.
\endproof

To proceed with the assumption ($\hat H_4$) when $f$ is not ``smoothifying'', recall that  a finite
measure $\varkappa$ is said to satisfy the \emph{Cramer's condition}  if
\begin{equation}
\sup_{|z|\geq \eps}\Big| \int_{\mathbb{R}} e^{iyz} \varkappa(dy)\Big|
<\varkappa(\mathbb{R})\quad \hbox{for  all }\, \eps>0 \label{Cr}
\end{equation}
(see,
for example, \cite{Lu79} or \cite{IL65}, chapter 3 \S3). Cramer's condition  means that  $\varkappa$  is in some sense regular. For instance, if $\varkappa$ has a non-trivial absolutely continuous part, then  (\ref{Cr}) follows from the Riemann-Lebesgue lemma, although, in general, a measure satisfying Cramer's condition should not be necessarily absolutely continuous (see Example \ref{ex1} below).

Note that (\ref{Cr}) leads to
$$
\Xi(\eps):=\sup_{|z|\geq \eps} \int_{\mathbb{R}} (1-\cos yz) \varkappa(dy)>0
\quad \hbox{for  all }\, \eps>0.
$$
In addition, assuming $\varkappa$ to have  finite second moment, we get  $\Xi(\eps)\sim c\eps^2$ as $\eps\to 0$ with some positive $c$, and thus
\begin{equation}
\Xi(\eps)=\sup_{|z|\geq \eps} \int_{\mathbb{R}} (1-\cos yz) \varkappa(dy)\geq c (\eps^2\wedge 1) \quad \hbox{for  all }\, \eps>0 \label{Cr1}
\end{equation}
 and some positive $c$.  Note that the function $\Theta(z, A)$ involved in ($\hat H_4$) is just the term under the supremum in (\ref{Cr1}), with $\varkappa$ equal to $\mu_f$ restricted to $A$. By the standing  assumptions  on $\mu$ and $f$, the measure $\mu_f$ restricted to $\Re\setminus(-q,q)$ has  finite second moment for any $q>0$. Therefore,  ($\hat H_4$) holds true,  provided that  for some $q>0$ the  restriction of  $\mu_f$  to  $[r,+\infty)$ satisfies the Cramer's condition.

\begin{lem}\label{l26}  Assume in addition to standing assumptions   on $\mu$ and $f$ that

($C$) for some $\rho>0$ the restriction of  $\mu$  to  $[\rho,+\infty)$ satisfies the  Cramer's condition.

Then  ($\hat H_4$) holds true for $q>0$ small enough.

\end{lem}

\proof  Take $r<\gamma  F\rho$ with $F=\mathrm{esssup}_{s\in I} f(s)$ and some $\gamma\in (0,1)$ . Then
$$\ba
\Theta(z,[r,+\infty))&=\iint_{(s,u):uf(s)\geq r}(1-\cos(uf(s)z))\mu(du)ds \\& \geq\int_{f(s)>\gamma F}\int_{u\geq \rho} (1-\cos(uf(s)z))\mu(du)ds \geq \left(\int_{f(s)>\gamma F}\,ds\right)\Xi_\rho(\gamma F\eps), \quad |z|\geq \eps,
\ea
$$
with $\Xi_\rho(\eps)=\sup_{|z|\geq \eps} \int_\rho^\infty (1-\cos uz) \mu(du)$. Since  $\Xi_\rho(\eps)$ satisfies (\ref{Cr1}) and the set $\{s: f(s)>\gamma F\}$ has positive Lebesgue measure, we obtain the required estimate for $\Theta(z,[q,+\infty))$.
\endproof

To summarise, let us formulate in the fixed time setting  the asymptotic results for the distribution densities of particular processes, listed in the Introduction.

\begin{cor}\label{c31} Let $Y$ be a L\'evy driven stochastic integral, specified below. Assume that the L\'evy measure of the noise  satisfies (\ref{mu0}),  (\ref{esint}), and ``tail'' conditions ($T_1$), ($T_2$).

 Then for every $t>0$ the distribution density $p_t$  exists, belongs to $C_b^\infty$, and satisfies
\begin{equation}\label{answ2}
p_t(x)\sim \frac{1}{\sqrt{2\pi K(t,x)}} e^{D(t,x)}, \quad
x \to \infty,
\end{equation}
with respective functions $K(t,x)$ and $D(t,x)$, in the following cases:

\begin{itemize}\item[\it (1)] $Y$ is the L\'evy process $Z$, $\mu$ satisfies ($N_1$) and ($C$);
\item[\it (2)] $Y$ is the non-stationary version of a L\'evy driven Ornstein-Uhlenbeck process, $\mu$ satisfies ($N_2$);
\item[\it (3)] $Y$ is the stationary version of a L\'evy driven Ornstein-Uhlenbeck process, $\mu$ satisfies ($N_3$) (in that case, $p_t(x)$, $K(t,x)$ and $D(t,x)$ actually don't depend on $t$);
\item[\it (4)] $Y$ is the fractional  L\'evy motion.
\end{itemize}
\end{cor}

\subsection{Examples}

In this section we give several examples that illustrate the
conditions on the measure $\mu$, introduced above.

The first two examples illustrate two typical situations
where ``tail'' conditions ($T_1$) and ($T_2$) hold.

\begin{exa}\label{ex2} Let $\mu$ be supported {in}  a bounded subset.
Denote by  $\sigma_+$ the minimal positive
constant $\sigma$ such that $\mu((\sigma,+\infty))=0$.  One can easily show
that for all  $\eps>0$ and  $k\geq 1$ one has
    \be\label{pm}
    M_k(\xi)\gg
    e^{(\sigma_+-\eps)\xi},\quad M_k(\xi)-\sigma^k_+
    \mu(\{\sigma_+\})e^{\sigma_+\xi}\ll e^{\sigma_+ \xi},
    \quad \xi\to +\infty.
    \ee
This relation yields both
($T_1$) and ($T_2$).  Indeed, for $\eps>0$ small enough one has
 $\gamma:=\frac{\sigma_{+}+\varepsilon}{2(\sigma_+ -\varepsilon)}\in (0,1)$ and
    $$
    M_4(\xi)\ll e^{(\sigma_++\varepsilon)\xi} = e^{2(\sigma_+-\varepsilon)\gamma \xi}\ll M_2^2(\gamma \xi),
    $$
which is ($T_1$). Similarly, for $\xi$ large enough
    $$
    \ln \left(\frac{M_4(\xi)}{M_2(\xi)}\vee 1 \right) +\ln \ln M_2(\xi)\leq  2\varepsilon \xi+\ln \xi,
        $$
which provides ($T_2$) because $\eps>0$ is arbitrary.

\end{exa}

\begin{exa}\label{ex3} Assume that  for $u$ large enough
\be\label{mu1}
{1\over Q(u)}e^{-b u^\beta}\leq \mu([u,+\infty))\leq Q(u)e^{-b u^\beta},
\ee
where $b>0$, $\beta>1$ are some constants,  and $Q$ is some polynomial.

For  $\sigma>0$ denote
$$
M_{k}^\sigma(\xi):=\int_{[\sigma,+\infty)}u^ke^{\xi u}\mu(du).
$$
Clearly,
$$
M_k^\sigma(\xi)\gg e^{A\xi},\quad  \xi \to +\infty
$$
for any $A>0$, and
$$
M_{k}(\xi)-M_{k}^\sigma(\xi)\ll e^{\sigma \xi},\quad  \xi \to +\infty.
$$
This means that, for any $\sigma>0$
\be\label{i300}
\begin{split}
M_{k}(\xi)&\sim \int_{[\sigma,+\infty)}u^ke^{\xi u}\mu(du)
\\&
=\sigma^ke^{\xi \sigma}\mu([\sigma, +\infty))+\int_{[\sigma,+\infty)}\Big[ku^{k-1}+\xi u^{k}\Big]e^{\xi u}\mu([u, +\infty))\, du.
\end{split}
\ee
For any $\sigma>0$, $m \in\mathbb{Z}$, we have
\be\label{i30}
\int_{\sigma}^{\infty}u^me^{\xi u}e^{-b u^\beta}\,du\sim c_1(\beta, b, m)\xi^{\frac{2m+2-\beta}{2(\beta-1)}}
e^{c_{2}(\beta)b^{1\over \beta-1}\xi^{\frac{\beta}{\beta-1}}},\quad
\xi\to+\infty,
\ee
where $c_2(\beta)=\beta^{-\frac{1}{\beta-1}}-\beta^{-\frac{\beta}{\beta-1}}$ (we have no need to specify the constant $c_1(\beta, b, m)$). One can prove (\ref{i30}) applying the Laplace method in a standard way; we omit the detailed calculations.

Take $\sigma$ large enough; then  (\ref{mu1}) holds true for $u\geq \sigma$. Then (\ref{i300}) and (\ref{i30}) yield for every $k\geq 1$,
\begin{equation}
{1\over Q_k(\xi)}e^{c_{2}(\beta)b^{1\over \beta-1}\xi^{\frac{\beta}{\beta-1}}}\leq  M_k(\xi)\leq Q_k(\xi) e^{c_{2}(\beta)b^{1\over \beta-1}\xi^{\frac{\beta}{\beta-1}}}\label{i3}
\end{equation}
for $\xi$ large enough, where $Q_k$ is some polynomial.

 By \eqref{i3},
    $$
    M_4(\xi) \leq Q_4(\xi) e^{c_{2}(\beta)b^{1\over \beta-1}\xi^{\frac{\beta}{\beta-1}}}\ll\left( {1\over Q_2(\gamma \xi)}e^{c_{2}(\beta)b^{1\over \beta-1}(\gamma \xi)^{\frac{\beta}{\beta-1}}}\right)^2\leq M_2^2(\gamma \xi),\quad  \xi \to +\infty,
    $$
as soon as $2\gamma^{\frac{\beta}{\beta-1}}>1$, which implies ($T_1$). Further, for $\xi$ large enough (\ref{i3}) gives
    \begin{align*}
    \ln &\left(\frac{M_4(\xi)}{M_2(\xi)}\vee 1 \right) +\ln \ln M_2(\xi) \\&
    \leq \ln\Big(Q_2(\xi)Q_4(\xi)\vee 1 \Big)+  \ln \Big(c_{2}(\beta)b^{1\over \beta-1}\ln Q_2(\xi)\Big)
     +{\beta\over \beta -1}\ln \xi
    \ll \xi,\quad \xi\to\infty,
    \end{align*}
which provides ($T_2$).

\end{exa}

The following example  illustrates condition ($C$) and the relations between the conditions ($N_1$) -- ($N_4$). All the measures in this example have bounded supports, therefore ``tail'' conditions ($T_1$) and ($T_2$) are satisfied.

\begin{exa}\label{ex1}
(a) Let $\mu=\sum_{n=1}^\infty n^{\rho} \delta_{n^{-1}}$, $\rho<1$; the assumption on $\rho$ provides  that $\mu$ is a L\'evy measure.
In the case $\rho\in (-1,1)$, the asymptotic behavior of the  integrals   $\int_{|u|\leq \eps}u^2\mu(du)$ is the same as for the $\alpha$-stable case with $\alpha=1+\rho$, i.e. is of a power type:
\be\label{orey}
\int_{|u|\leq \eps}u^2\mu(du)\asymp \eps^{2-\alpha}, \quad \eps\to 0
\ee
(cf. \cite{O68}, \cite{Pi96}, \cite{IK05}). Therefore, conditions ($N_1$) -- ($N_4$)  holds true. The analogy with the $\alpha$-stable case is not complete: condition $(C)$ does not hold, because for every  $r>0$ the restriction of $\mu$ to $[r,+\infty)$  has finite number of atoms. Therefore,  statements (2) -- (4) in Corollary \ref{c31} hold true, but one can not claim (\ref{answ2}) for the L\'evy process $Z$ itself. Statement (1) of Corollary \ref{c31} becomes applicable when
$\mu$ is replaced by $\mu+\varkappa,$ where $\varkappa$ is a measure with a bounded support, satisfying Cramer's condition. For instance, either $\kap$ may  be absolutely continuous  (and then Cramer's condition is provided by the Riemann-Lebesgue lemma), or $\kap$ may be equal to the Cantor measure on $[0,1]$ (and then Cramer's condition is verified by  straightforward calculations).

When $\rho=-1$, ($N_1$) and ($N_2$) fail, but ($N_3$) and ($N_4$) hold true. When $\rho<-1$, only ($N_4$) hold true, while ($N_1$) -- ($N_3$) fail. It is clear that, in the latter case,  the laws of the L\'evy process $Z$ and of the non-stationary version of a L\'evy driven Ornstein-Uhlenbeck process contain non-trivial discrete components. Therefore one definitely can not expect any asymptotic relation like (\ref{answ2}) to hold for these processes. On the other hand, (\ref{answ2}) holds true for the fractional L\'evy motion $Z^H$ with $H\in (1/2,1)$. This well illustrates the ``smoothifying'' role of the kernel $f$.

(b) $\nu = \sum_{n=1}^\infty n \delta_{({n!})^{-1}}$. Then condition $(N_1)$ fails,  while
($N_2$) -- ($N_4$) hold true, see  \cite{Ku06}, Example 2.3 or
\cite{BK08}, Example 1.   In these examples, it is shown  that the law of
$Z_t$ is singular for {all}  $t>0$. Thus,  the asymptotic relation (\ref{answ2}) clearly can not be valid for the L\'evy process $Z$ itself. In this case, the L\'evy measure provides some ``hidden smoothness'' in the sense that the law of the L\'evy process $Z$ is singular, but the distributions of the respective (both non-stationary and stationary) L\'evy driven Ornstein-Uhlenbeck processes  and   fractional L\'evy motion possess $C^\infty_b$ distribution densities which, moreover, admit asymptotical description  (\ref{answ2}).

\end{exa}

In the last example, in the case of a L\'evy process, we compare our conditions with those introduced  in \cite{CKK08} and \cite{CKK10}.

\begin{exa}\label{ex4}

In the paper \cite{CKK08}  the authors give the transition density estimates for  a symmetric
$\alpha$-stable--like process whose jump intensity kernel $J(x,y)$  is of
the form
    $$
    J(x,y)=\frac{c(x,y)}{|x-y|^{n+\alpha}}1_{|x-y|\leq1},
    $$
where $c(x, y)$ is a symmetric Borel measurable function
on $\mathbb{R}^n\times \mathbb{R}^n$,  bounded from above and
below by  two positive constants. When $c(x,y)\equiv c(x-y)$, that is, $J(x,y)=J(x-y)$, this process is a L\'evy one with the L\'evy measure $\mu(dx)=J(x)dx$. We check in {the} one-dimensional case that such a L\'evy measure satisfies the conditions imposed above.

Since the  L\'evy measure $\mu$   has bounded support, the exponential integrability condition is satisfied and, moreover, conditions ($T_1$) and ($T_2$) hold true (see Example~\ref{ex2}). By the Riemann-Lebesgue lemma, the absolute continuity of $\mu$ implies condition ($C$).  Finally,
 (\ref{orey}) holds true, which provides the Kallenberg condition \eqref{Kal} for the measure $\mu$. Since $\mu$ is assumed to be symmetric, this yields ($N_1$).

The paper  \cite{CKK10}  is devoted to the estimates of the transition density  of a Markov process whose jump intensity  $J(x,y)$ satisfies
    \begin{equation}
    \frac{c_1}{|x-y|^n\phi(c_2|x-y|)} \leq J(x,y)\leq  \frac{c_3}{|x-y|^n\phi(c_4|x-y|)}, \quad x,y\in \mathbb{R}^n\times\mathbb{R}^n, \quad x\neq y,\label{j2}
    \end{equation}
for some $c_i$, $i=1,2,3,4$, where $\phi: [0,\infty)\to [0,\infty)$  is of the form $\phi(r)=\phi_1(r)\psi(r)$, $r>0$, and

i) $\psi$ is increasing  on $[0,\infty)$, $\psi(r)=1$ for $0<r\leq 1$,  and for some  $0<\gamma_1\leq \gamma_2$, $\beta>0$,
    \begin{equation}
    c_1e^{\gamma_1 r^\beta}\leq \psi(r)\leq c_2 e^{\gamma_2 r^\beta},\quad 1<r<\infty;\label{psi1}
    \end{equation}

ii) $\phi_1$ is strictly increasing on $[0,\infty)$ with $\phi_1(0)=0$, $\phi_1(1)=1$, and, in particular,  satisfies for $c_2>c_1>0$, $c_3>0$, $0<\beta_1\leq \beta_2<2$, the inequality
    \begin{equation}
    c_1\left(\frac{R}{r}\right)^{\beta_1} \leq \frac{\phi_1(R)}{\phi_1(r)}\leq c_2\left(\frac{R}{r}\right)^{\beta_2}  \quad \text{for every $0<r<R<\infty$}.\label{phi1}
    \end{equation}

Again, let $n=1$ and $J(x,y)=J(x-y)$, where $J$ is the density of the L\'evy measure $\mu$.  To achieve the exponential integrability (\ref{esint}) we need to assume   $\beta>1$ in \eqref{psi1}. Assuming additionally that $\gamma_1=\gamma_2$,   one has   ($T_1$) and ($T_2$) (see Example~\ref{ex3}). By
  \eqref{phi1}, the L\'evy measure $\mu$ satisfies the lower bound in (\ref{orey}) with $\alpha=\beta_1$; that is,
    $$
  \int_{|u|\leq \eps}u^2\mu(du)\geq c\eps^{2-\beta_1}
  $$
  with some positive $c$ and $\eps>0$ small enough, which implies ($N_1$). Finally,   condition ($C$) holds true by the  absolute continuity of $\mu$.

Since ($N_1$), ($C$), ($T_1$), and ($T_2$) hold true, by statement (1) in Corollary \ref{c41} and Corollary \ref{c51} below, {the} transition probability density of the Levy process $Z$ satisfies (\ref{answ41}) and either (\ref{answ711}) (in the  ``truncated" case \cite{CKK08}) or (\ref{answ712})  (in the case treated in \cite{CKK10}). Let us compare these relations with the estimates for the transition probability density of a symmetric jump process from \cite{CKK08} and \cite{CKK10}.

For $t\geq t_0$, these estimates are given  in the form
\begin{equation}
C_1 g_t(C_2|x-y|)\leq p(t,x,y)\leq C_3 g_t(C_4|x-y|), \label{pg}
\end{equation}
where $C_1, \dots, C_4$ are some positive constants, and
\begin{equation}\label{g}
g_t(x)=\exp\left(-|x|\ln^\delta
\frac{|x|}{t}\right)\vee \left({1\over
t^{d/2}}\exp\left(-{|x|^2\over t}\right)\right)
\end{equation} with
$\delta=1$ in the ``truncated" case \cite{CKK08}, and
$\delta=\frac{\beta}{\beta-1}$ in the case treated in \cite{CKK10}
($d$ is the dimension of the space;  in the current paper $d=1$).

For a Levy process,  (\ref{pg}) with  $p(t,x,y)=p_t(y-x)$ is
closely comparable  with (\ref{answ41}) and (\ref{answ711}),
(\ref{answ712}).  When ${|x-y|\over t}$ is large, (\ref{answ711}),
(\ref{answ712}) directly provide (\ref{pg}) with $g_t$ replaced by
$$
\exp\left(-|x|\ln^\delta \frac{|x|}{t}\right).
$$
On the other hand, one can show easily that on every bounded set
the function  $\Kf_Z$ is bounded and bounded away from $0$, and
the function $\Df_Z$ satisfies
$$
-d_1x^2\leq \Df_Z(x)\leq -d_2x^2
$$
with positive constants $d_1, d_2$. Thus,   when ${|x-y|\over t}$
is  bounded, (\ref{answ41}) provides  (\ref{pg}) with $g_t$
replaced by
$$
{1\over  t^{1/2}}\exp\left(-{|x|^2\over t}\right).
$$
Note that (\ref{answ711}) and  (\ref{answ712}) are somewhat more
precise than (\ref{pg}): by choosing $x\over t$  large enough, one
can make the constants $c_1,c_2$ therein to  be arbitrarily close
to a given constant $c_*$,  while in (\ref{pg}) respective
constants $C_2$ and $C_4$ are different and fixed.

Although having a non-trivial intersection,  the classes  of
processes, treated in our case and in  \cite{CKK08} and
\cite{CKK10}, are substantially different. Our approach, based on
the Fourier transform technique, is not applicable to the class of
symmetric jump processes from \cite{CKK08} and  \cite{CKK10} in
the whole generality. On the other hand, this approach  is
applicable to particularly interesting processes which can not be
studied by the technique of \cite{CKK08}, \cite{CKK10}, including
non-symmetric Markov processes (like the L\'evy driven
Ornstein-Uhlenbeck process) and non-Markov processes (like the
fractional L\'evy motion).

\end{exa}

\section{Explicit conditions: time-dependent setting}\label{s4}

Our further aim is to consider conditions of Theorem \ref{td1} in the general, i.e.  time-dependent, setting.
To make the exposition reasonably short, we address this problem in a particular case of the self-similar kernel $f$; that is, we assume that
\be\label{self}
f(t,s)=\chi(t)f\left( {s\over\theta (t)} \right), \quad t\in \TT, \quad s\in I
\ee
with some functions $f:\Re \to \Re$ and $\chi,\theta:\TT\to (0,+\infty)$. Assumption (\ref{self}) is satisfied for particularly interesting  processes like the L\'evy process $Z$ and the fractional L\'evy motion $Z^H$.
In these cases we have, respectively,
\be\label{Z}
 f(s)=\1_{[0,1]}(s), \quad \chi(t)=1, \quad \theta(t)=t;
\ee
\be\label{ZH}
f(s)={1\over \Gamma(H+1/2)}\left[(1-s)_+^{H-1/2}-(-s)^{H-1/2}_+\right], \quad \chi(t)=t^{H-1/2}, \quad \theta(t)=t.
\ee

For the function $f(s)$ we keep our  standard standing
assumptions:  it is bounded and satisfies (\ref{fint}),
(\ref{fnon}). For the L\'evy measure
 $\mu$ we assume   (\ref{mu0}) and (\ref{esint}) to hold true, as before.

 Denote, similarly to Section \ref{s2},
  $$
\Theta(z, A)=\iint_{\{(s,u)\in \Re\times \Re:\,  f(s) u\in A\}} (1-\cos(f(s) z u))\mu(du)ds, \quad z\in \Re,
$$
$$
\Psi(z)=\int_\Re\int_{\Re}\left(e^{-izf(s)u}-1+izf(s)u\right)\,\mu(du)ds,  \quad  z\in \CC,
$$
$$
H(y,z)=iyz+\Psi(z), \quad \Mf_{k}(\zeta)={\prt^k\over \prt \zeta^k}\Psi(i\zeta), \quad
k\geq 1, \quad y\in \Re, \quad \zeta\in\mathbb{R}.
$$
Denote by $\zeta(y)\in \Re$ the unique solution  to the equation
\be\label{eq1}
{\prt\over \prt \zeta}H(y,i\zeta)=0,
\ee
and put
\be\label{kd}
\Df(y)=H(y,i\zeta(y)),\quad
\Kf(y)=\Mf_{2}(\zeta(y))={\prt^2\over \prt
\zeta^2}H(y,i\zeta)\Big|_{\zeta=\zeta(y)}.
\ee

Denote $\tau(t)=\chi(t)\theta(t)$.  Further in this section we
assume $\theta$ and $\chi$ to be bounded on every segment
$[a,b]\subset (0,+\infty)$, and to be bounded away from $0$ on the
whole $\TT$.  Clearly, the functions $\theta, \chi$ in (\ref{Z})
 and in (\ref{ZH}) with $H>1/2$ satisfy these
assumptions. In addition, we assume  that  \be\label{chitheta}
\theta(t) \to +\infty, \quad \ln\Big((\ln \chi(t))\vee 1\Big)\ll \ln
\theta (t), \quad t\to +\infty; \ee
 in the cases (\ref{Z}) and (\ref{ZH}) this assumption holds true.

\begin{theorem}\label{t41}

Assume that the measure $\mu$  satisfies   ($T_1$) and  ($T_2$). Assume also {that} $\mu$  {satisfies} one of the conditions ($N_i$) and, respectively, $f$ {satisfies} one of the assumptions ($F_i$), $i=1,\dots,4$. In the case $i=1$, assume additionally  that $\mu$  satisfies condition ($C$).

Then for every  $t>0$ the law of $Y_t$ has a distribution density $p_t\in C_b^\infty$, and for every $t_0>0$
\begin{equation}\label{answ3}
p_t(x)\sim {1\over \tau(t)}\sqrt{\frac{\theta(t)}{2\pi \Kf({x/ \tau(t)})}} e^{\theta(t)\Df\left({x/\tau(t)}\right)}, \quad
t+x\to \infty,\quad (t,x)\in [t_0, +\infty)\times \ax.
\end{equation}

\end{theorem}

\begin{rem} {The} expression {on} the the right hand side of (\ref{answ3}) is self-similar in the sense that the variable $x$, rescaled by $\tau(t)$, is involved in this expression only  as an argument of given functions $\Kf$ and $\Df$. Note that the L\'evy measure $\mu$ is not assumed to have a self-similarity property, and therefore, in general, the family of distributions $Y_t$, $t>0$ is not self-similar. Thus,  although the assumption (\ref{self}) on the kernel $f$ itself does not provide self-similarity for the \emph{distribution densities} of $Y_t$, $t>0$, it is powerful enough to provide self-similarity for the  \emph{asymptotic relation} for these densities.
\end{rem}

\begin{proof} The  relations below follow easily from the self-similarity assumption (\ref{self}):
\be\label{e2}
H(t,x,z)=\theta(t)H\left({x\over \tau(t)}, \chi(t)z\right), \quad \Mf_k(t,\xi)=\chi^k(t)\theta(t)\Mf_k(\chi(t)\xi), \quad k\geq 1.
\ee
By the first relation in (\ref{e2}),  we can rewrite the relation  (\ref{eq}), which determines $\xi=\xi(t,x)$,  as
$$
\chi(t)\theta(t) {\prt\over \prt\zeta} H\left({x\over \tau(t)}, i\zeta\right)\Big|_{\zeta=\chi(t)\xi}=0.
$$
This means that $\chi(t)\xi$ solves (\ref{eq1}) with $y=x/\tau(t)$, and therefore
$$
\xi(t,x)=\chi^{-1}(t)\zeta\left({x\over \tau(t)}\right).
$$
Combined with (\ref{e2}), this relation gives
$$
D(t,x)=\theta(t)\Df\left({x\over \tau(t)}\right), \quad K(t,x)=\chi^2(t)\theta(t)\Kf\left({x\over \tau(t)}\right)={\tau^2(t)\over \theta(t)}\Kf\left({x\over \tau(t)}\right).
$$
Thus (\ref{answ3}) would follow from (\ref{answ1}) with $\Af=[t_0, +\infty)\times\Re^+$, provided that conditions ($H_1$) -- ($H_4$) are verified.

In Section 3  we proved that under assumptions imposed on  the L\'evy measure
$\mu$ and the function $f(s)$,  conditions ($\hat H_1$),  ($\hat H_2$), ($\hat H_3^s$), and  ($\hat H_4$) hold true.
Now we show that these conditions yield ($H_1$) -- ($H_4$) with
  $\mathcal{T}=[t_0,+\infty)$, $\mathcal{B}=[t_0,+\infty)\times \Re^+$, and {with the} function $\theta(t)$ replaced by  $\vartheta\theta(t)$ (the constant $\vartheta$ comes from ($\hat H_4$)).

The second relation in (\ref{e2}) gives
\be\label{e3}
{\Mf_4(t,\xi)\over \Mf_2^2(t,\xi)}={1\over
\theta(t)}{\Mf_4(\chi(t)\xi)\over \Mf_2^2(\chi(t)\xi)}.
\ee
Observe that, under our  assumptions on $\theta$ and  $\chi$,
\be\label{e5}
 t+\xi\to \infty\quad \text{implies } \quad \theta(t)\to
+\infty\quad \hbox{or}\quad \chi(t)\xi\to +\infty.
\ee
 Therefore,
($H_1$) follows from ($\hat H_1$) and (\ref{e3}).

By the second relation in (\ref{e2}),
\be\label{e4}
{\Mf_4(t,\xi)\over \Mf_2(t,\xi)}=\chi^2(t){\Mf_4(\chi(t)\xi)\over
\Mf_2(\chi(t)\xi)},
\ee
which together with ($\hat H_2$) and
(\ref{e5}) gives
$$
\ln\left(\Big(\chi^{-2}(t)\frac{\Mf_{4}(t,\xi)}{\Mf_{2}(t,\xi)}\Big)
\vee 1\right)\ll \ln \theta(t)+\chi(t)\xi, \quad  t+\xi\to +\infty.
$$
Similarly,
$$\ba
\ln\left(\Big( \ln
\Mf_{2}(t,\xi)\Big)\vee 1\right)&=\ln\left(\ln\Big(\chi^2(t)\theta(t)\Mf_{2}(\chi(t)\xi)\Big)\vee 1\right)\\
& = \ln\left(\Big(\ln\chi^2(t)+\ln\theta(t)+\ln
\Mf_{2}(\chi(t)\xi)\Big)\vee 1\right).\ea
$$
By ($\hat H_2$),  (\ref{e5}) and (\ref{chitheta}) one has
$$
\ln\left(\Big( \ln
\Mf_{2}(t,\xi)\Big)\vee1\right)\ll \ln \theta(t)+\chi(t)\xi, \quad  t+\xi\to +\infty.
$$
This completes the proof of ($H_2$).

By ($\hat H_3^s$), for every $\varkappa>0$ there exists $Q>0$ such that
$$
\Theta(z, \Re^+) \geq \varkappa \ln|z|, \quad |z|\geq Q.
$$
By  the self-similarity assumption (\ref{self}), we have
$$
\Theta(t,z,A)={\theta(t)}\Theta\left(\chi(t)z,{1\over \chi(t)} A\right).
$$
Denote $\theta_*=\inf_{t}\theta(t), \chi_*=\inf_{t}\chi(t)$. Then  taking  $\varkappa=\theta_*^{-1}(1+\delta)$ and $R=\chi_*^{-1}Q$, we obtain ($H_3$).

Finally,   by ($\hat H_4$) we have
$$\ba
\inf_{t\in \mathcal{T}, \, |z|>\eps}\Theta(t,z,[ q\chi(t)< +\infty) )&=\theta(t)\inf_{t\in \mathcal{T},\,  |z|>\eps}
\Theta(\chi(t) z, [q< +\infty) )\\&=
\theta(t)\inf_{t\in \mathcal{T}, \, |z'|>\chi(t)\eps}\Theta(z', [q< +\infty) )\geq \vartheta \theta(t)\Big((\chi(t)\eps)^2\wedge 1\Big).\ea
$$
Thus,  ($H_4$) holds true with $r=q$ and $\theta(t)$ replaced by $\vartheta\theta(t)$. Clearly, such a change of the function $\theta(t)$ does not spoil conditions ($H_1$) -- ($H_3$) proved above.
\end{proof}

\begin{cor}\label{c41}  Assume the L\'evy measure of the noise satisfy (\ref{mu0}),  (\ref{esint}), and ``tail'' conditions ($T_1$), ($T_2$). Then

\begin{itemize}\item[\it (1)] For the L\'evy process $Z$, assuming additionally $\mu$ to satisfy ($N_1$) and ($C$), one has
(\ref{answ41}).
\item[\it (2)] For the fractional  L\'evy motion $Z^H$, one has (\ref{answ42}).
\end{itemize}
\end{cor}

\section{Explicit asymptotic expressions as $x\to+\infty$}\label{s5}

Theorem \ref{td1}, Corollary \ref{c31}, and Theorem \ref{t41}    describe the asymptotic
behaviour  of
a distribution density precisely, but in an implicit form:  functions $K(t,x)$, $D(t,x)$, $\Kf(x)$, $\Df(x)$, involved in (\ref{answ1}), (\ref{answ2}) and  (\ref{answ3}),  are defined in  terms of  the solutions to equations (\ref{eq}) or (\ref{eq1}). In this section we
study the asymptotic behavior of these functions as $x\to+\infty$, and
 deduce  explicit  asymptotic
expressions for the distribution densities.

In what follows, we mainly discuss the behavior of the functions
$\Kf(x)$ and $\Df(x)$ under additional  assumptions on the L\'evy
measure $\mu$; without any essential change of the argument,
similar results can be obtained for the functions $K(t,x)$,
$D(t,x)$ with a fixed variable $t$.  To simplify the argument,
we assume in the sequel $f$ to be non-negative. This assumption
is satisfied, for instance, for the L\'evy process and the
fractional L\'evy motion (respective functions $f$ are given in
(\ref{Z}) and (\ref{ZH})). To shorten the exposition, we restrict
ourselves  to the cases where the  L\'evy measure $\mu$ is either
``truncated'' (i.e. supported in a bounded set, see Example
\ref{ex2}) or ``exponentially damped'' (i.e. having its ``tails''
satisfying  (\ref{mu1}), see Example \ref{ex3}).

We keep the notation introduced in Example \ref{ex2}, Example \ref{ex3}, and  Section \ref{s4}; in particular, $F=\mathrm{essup}\, f(s)$, and $\sigma_+$ is the extreme right point  of the support of $\mu$.

\begin{theorem}\label{t51} Assume  the kernel $f(t,s)$ to be of the form (\ref{self}) with $\theta$ and $\chi$ satisfying (\ref{chitheta}). Assume that the measure $\mu$  satisfies   one of the conditions ($N_i$), and the respective function  $f(s)$ in (\ref{self})  satisfies  one of the assumptions ($F_i$), $i=1,\dots,4$. In the case $i=1$, assume additionally $\mu$ to satisfy condition ($C$).

 1.  If $\mu$ is truncated, then   for any constants $c_1>1/(\sigma_+F)$ and $c_2<1/(\sigma_+F)$ there exists $y=y(c_1,c_2)$ such that, for  $x/\tau(t)>y$,
 \be\label{answ51}
 \exp\left(-c_1\left(x\over \chi(t)\right)\ln \left(x\over \tau(t)\right)\right)\leq p_t(x)\leq \exp\left(-c_2\left(x\over \chi(t)\right)\ln \left(x\over \tau(t)\right)\right).
\ee

 2.   If $\mu$ is exponentially damped, then for any constants
 $$
    c_2<\left(\beta^{-\frac{1}{\beta-1}}-\beta^{-\frac{\beta}{\beta-1}}\right)^{-{\beta-1\over \beta}}b^{-{1\over \beta}}F^{-1}<c_1
 $$
 there exists $y=y(c_1,c_2)$ such that, for  $x/\tau(t)>y$,
\be\label{answ52}
\exp\left(-c_1\left(x\over \chi(t)\right)\ln^{\beta-1\over \beta} \left(x\over \tau(t)\right)\right)\leq p_t(x)\leq \exp\left(-c_2\left(x\over \chi(t)\right)\ln^{\beta-1\over \beta} \left(x\over \tau(t)\right)\right).
\ee

\end{theorem}

\proof We consider in {detail} the case of a truncated L\'evy measure, and then outline the changes in the proof that should be made in the case of an exponentially damped L\'evy measure.

Denote
$$
\Mf_0(\zeta)=\Psi(i\zeta)=\int_\Re\int_{\Re}\left(e^{\zeta f(s)u}-1-\zeta f(s)u\right)\,\mu(du)ds.
$$
Similarly to (\ref{m2}), one has
$$
\Mf_0(\zeta)=\int_{\Re}M_0(f(s)\zeta)\, ds, \quad M_0(\xi):=\int_{\mathbb{R}} \big(e^{\xi u}-1-\xi u
\big) \mu(du).
$$

To  describe the asymptotic behavior of $\Kf$, $\Df$, we need to analyze the behavior of $\Mf_k$, $k=0,1,2$.  For this, we analyze first the behavior of $M_k$, $k=0,1,2$.

One can easily see that (\ref{pm}) holds true for $k=0$ as well. From (\ref{pm}) we have for any $k\geq 0$
\be\label{73}
M_k(\xi)\sim \sigma_{+}^k M_0(\xi), \quad \xi \to
+\infty.
\ee
Moreover, the first relation in (\ref{pm}) provides that for every $\eps>0$
\be\label{730}
\Mf_k(\zeta)\sim \int_{f(s)\geq F-\eps }f^k(s)M_k(f(s)\zeta)\, ds, \quad \zeta\to +\infty,
\ee
(recall that we assume $f$ to be non-negative), which together with (\ref{73}) yields
\be\label{731}
\Mf_k(\zeta)\sim F^k\sigma_{+}^k\Mf_0(\zeta), \quad \zeta \to
+\infty.
\ee

Recall that $\Kf(y)=\Mf_2(\zeta(y))$, and
$$
\Df(y)=-y\zeta(y)+\Mf_0(\zeta(y)).
$$
The function $\zeta(y)$ is defined by the equation  $\Mf_1(\zeta(y))=y$ and, under our fixed  assumption (\ref{mu0}), we have  $\zeta(y)\to +\infty$ as $y\to +\infty$. Hence, by (\ref{731}),
\be\label{kf}
\Kf(y)\sim F\sigma_+ y, \quad \Mf_0(\zeta(y))\sim (1/F\sigma_+) y, \quad y\to
+\infty.
\ee
The second relation in the above formula yields
\be\label{df}
\Df(y)\sim -y\zeta(y),   \quad y\to
+\infty.
\ee

Similarly to (\ref{731}), one can deduce from (\ref{pm}) that for any $\eps>0$
$$
 e^{(\sigma_+ F-\eps)\zeta}\ll \Mf_1(\zeta)\ll e^{(\sigma_+ F+\eps)\zeta}, \quad \zeta \to +\infty,
$$
and consequently
\be\label{zeta0}
\zeta(y)\sim {1\over \sigma_{+}F}\ln y, \quad y \to +\infty.
\ee

 Let us prove the lower bound in (\ref{answ51}), the proof of the upper bound is similar and omitted. It follows from (\ref{df}) and (\ref{zeta0}) that for any $c> 1/(\sigma_+F)$ we have  for  $x/\tau(t)$ large enough
\be\label{51}
e^{\theta(t)\Df(x/\tau(t))}\geq \exp\left(-c\left(\theta(t) x\over \tau(t)\right)\ln \left(x\over \tau(t)\right)\right)=\exp\left(-c\left(x\over \chi(t)\right)\ln \left(x\over \tau(t)\right)\right),
\ee
(recall that $\tau(t)=\theta(t)\chi(t)$).

 Since $\mu$ is supported in a bounded set,   it satisfies ``tail'' conditions ($T_1$), ($T_2$) (see  Example \ref{ex2}), and thus we can apply  Theorem \ref{t41}. By Theorem \ref{t41} and (\ref{51}), to prove the first inequality in (\ref{answ51}) it is enough to take $c\in \Big(1/(\sigma_+F), c_1\Big)$ and prove that for $x/\tau(t)$ large enough,
\be\label{52}
{1\over \tau(t)}\sqrt{\frac{\theta(t)}{2\pi \Kf({x/ \tau(t)})}}\geq \exp\left((c-c_1)\left(x\over \chi(t)\right)\ln \left(x\over \tau(t)\right)\right).
\ee
By (\ref{kf}) and (\ref{chitheta}), we have  for any $q>0$,
\be\label{53}
{1\over \tau(t)}\sqrt{\frac{\theta(t)}{2\pi \Kf({x/ \tau(t)})}} \geq  {1\over \theta(t)}e^{-q x/ \tau(t)}
\ee
for $x/\tau(t)$ large enough. On the other hand, for a fixed $y>1$ and $x/\tau(t)\geq y$,
$$
\exp\left((c-c_1)\left(x\over \chi(t)\right)\ln \left(x\over
\tau(t)\right)\right) \leq  e^{-q\theta(t)(x/\tau(t))}\quad
\hbox{with}\quad q=(c_1-c_2)\ln y>0.
$$
Hence (\ref{52}) follows from the inequality
$$
{1\over a}e^{-qb}=e^{-\ln a-qb}\geq e^{-q(a+b)}\geq e^{-q ab},
$$
valid for $a,b$ large enough.

Let us discuss briefly the changes that should be made when the measure $\mu$ satisfies (\ref{mu1}). Clearly, for any $\sigma>0$ and  $k>j$,
$$
M_k^\sigma(\xi) \geq \sigma^{k-j} M_j^\sigma(\xi), \quad \xi\geq 0
$$
(see the notation in Example \ref{ex3}). This means that instead of (\ref{73}) and (\ref{731}) we have now
\be\label{732}
M_k(\xi)\gg M_j(\xi), \quad \xi \to
+\infty,
\ee
\be\label{733}
\Mf_k(\zeta)\gg \Mf_j(\zeta), \quad \zeta \to
+\infty
\ee
for any $k>j$. The latter relation with $k=1, j=0$ yields (\ref{df}). From (\ref{i3}) and (\ref{730}) it follows that for every $\eps>0$ for $\zeta$ large enough $$
e^{(C_*-\eps)\zeta^{\frac{\beta}{\beta-1}}}\leq  \Mf_1(\zeta)\leq e^{(C_*+\eps)\zeta^{\frac{\beta}{\beta-1}}},
$$
where
$C_*=\left(\beta^{-\frac{1}{\beta-1}}-\beta^{-\frac{\beta}{\beta-1}}\right)b^{1\over
\beta-1}F^{\frac{\beta}{\beta-1}}$. Consequently, \be\label{54}
\zeta(y)\sim \left({1\over C_*}\ln y\right)^{\beta-1\over \beta},
\quad \Df(y)\sim -y\left({1\over C_*}\ln y\right)^{\beta-1\over
\beta}, \quad y\to +\infty, \ee which  means that the  analogue
of (\ref{answ52}), with $e^{\theta(t)\Df(x/\tau(t))}$ instead of
$p_t(x)$, holds true,  and the only thing we need to verify is
that the term
$$
{1\over \tau(t)}\sqrt{\frac{\theta(t)}{2\pi \Kf({x/ \tau(t)})}}
$$
is  negligible.   Note that this term is bounded:
$$
\sup_t {1\over \tau(t)}\sqrt{\frac{\theta(t)}{2\pi }}=\sup_t{1\over \sqrt{2\pi \chi^2(t)\theta(t)}}<+\infty
$$
because $\theta$ and $\chi$ are assumed to be separated from $0$, and by  (\ref{733})
$$
\Kf({x/ \tau(t)})=\Mf_2(\zeta({x/ \tau(t)}))\gg \Mf_1(\zeta({x/ \tau(t)}))={x/ \tau(t)}, \quad {x/ \tau(t)}\to +\infty.
$$
 This observation provides the upper bound in (\ref{answ52}).

On the other hand, it follows from (\ref{i3}) that
$$
\ln {M_2(\xi)\over M_1(\xi)}\ll \xi, \quad \xi \to +\infty
$$
(cf. (\ref{732})). Similarly to the proof of Lemma \ref{l23}, one can deduce from this relation that
$$
\ln {\Mf_2(\zeta)\over \Mf_1(\zeta)}\ll \zeta, \quad \zeta \to +\infty,
$$
and consequently
$$\ln \Kf(y)\ll \zeta(y)+\ln y, \quad y\to +\infty.
$$
Together with (\ref{chitheta}) and the first relation  in (\ref{54}), this implies (\ref{53}). Repeating the argument after (\ref{53}), we obtain the lower bound in (\ref{answ52}).
\endproof

For the L\'evy process $Z$ and the fractional L\'evy motion  $Z^H$,  Theorem \ref{t51} gives the following. Denote
\be\label{cstar}
c_*= 1/\sigma_+
\ee
in the case of the  truncated L\'evy measure $\mu$, and
\be\label{cstar1}
c_*= \left(\beta^{-\frac{1}{\beta-1}}-\beta^{-\frac{\beta}{\beta-1}}\right)^{-{\beta-1\over \beta}}b^{-{1\over \beta}}
\ee
in the case of the exponentially damped L\'evy measure $\mu$.

\begin{cor}\label{c51}  Assume that the L\'evy measure satisfies  ($N_1$) and ($C$), then for the distribution density of the L\'evy process $Z$ the following estimates hold.

1.  If $\mu$ is truncated, then   for any constants $c_1>c_*$ and $c_2<c_*$ there exists $y=y(c_1,c_2)$, such that for  $x/t>y$ (\ref{answ711}) holds true.

 2.   If $\mu$ is exponentially damped, then for any constants
 $c_1>c_*$ and $c_2<c_*$
 there exists $y=y(c_1,c_2)$, such that for  $x/t>y$, (\ref{answ712}) holds true.
\end{cor}

\begin{cor}\label{c52}  Assume that the L\'evy measure satisfies  (\ref{mu0}), then for the distribution density of the fractional L\'evy motion  $Z^H$ the following estimates hold.

1.  If $\mu$ is truncated, then   for any constants $c_1>c_*$ and $c_2<c_*$ there exists $y=y(c_1,c_2)$, such that for  $x/t^{H+1/2}>y$ (\ref{answ721}) holds true.

 2.   If $\mu$ is exponentially damped, then for any constants
$c_1>c_*$ and $c_2<c_*$
 there exists $y=y(c_1,c_2)$, such that for  $x/t^{H+1/2}>y$ (\ref{answ722}) holds true.
\end{cor}

We have mentioned in the beginning of the section that  for fixed $t$ the functions $K(t,x)$, $D(t,x)$ can be analyzed in the same way as $\Kf(x)$ and $\Df(x)$. Respectively, the analogue of Theorem \ref{t51} can be proved for the density $p_t(x)$  with  fixed $t$ without the self-similarity assumption (\ref{self}).  Let us formulate one statement of such a kind.

Consider the stationary version $X$ of a L\'evy driven
Ornstein-Uhlenbeck process, and assume  that $\mu$  satisfies
($N_3$). Then  the  distribution of $X_t$, in fact,
does not depend on $t$,  and  by Proposition \ref{p21} has a
$C^\infty$ distribution density $p$, which we call the
\emph{invariant distribution density} of the  respective
Ornstein-Uhlenbeck process. Moreover, assuming the ``tail''
conditions ($T_1$) and  ($T_2$) to hold, we have  by  Corollary
\ref{c31}  the asymptotic relation for this density, which after
trivial transformations can be written in the form
\begin{equation}\label{answ6}
p(x)\sim \frac{1}{\sqrt{2\pi \Kf(x)}} e^{\Df(x)}, \quad
x \to +\infty,
\end{equation}
where $\Kf$, $\Df$ are defined by (\ref{kd}) with $f(s)=e^{\gamma s}\1_{s\leq 0}$. Similarly to Theorem \ref{t51}, one can deduce from (\ref{answ6}) the following statement (the proof is omitted).

\begin{prop}\label{p51}  Assume that the L\'evy measure $\mu$ satisfies  (\ref{esint}) and ($N_3$). Then for the invariant distribution density of the  Ornstein-Uhlenbeck process (\ref{sde}) the following estimates hold.

1.  If $\mu$ is truncated, then   for any constants $c_1>c_*$ and $c_2<c_*$ there exists $y=y(c_1,c_2)$,  such that for  $x>y$,
$$
\exp\left(-c_1 x\ln x\right)\leq p(x)\leq \exp\left(-c_2 x\ln x\right).
$$

 2.   If $\mu$ is exponentially damped, then for any constants
 $c_1>c_*$ and $c_2<c_*$
 there exists $y=y(c_1,c_2)$ such that for  $x>y$,
$$
\exp\left(-c_1x\ln^{\beta-1\over \beta} x\right)\leq p(x)\leq \exp\left(-c_2x\ln^{\beta-1\over \beta} x \right).
$$

Here $c_*$ depends on $\mu$ only, and is defined, respectively,  in (\ref{cstar}) or (\ref{cstar1}).
\end{prop}

As we mentioned in the Introduction, there is a particular theoretical interest in studying  the ratio (\ref{c1}) of the values of the invariant distribution density. One can see that the statement of Proposition \ref{p51}  is not strong enough to provide an exact estimate for the ratio  (\ref{c1}) because of different constants $c_1$ and $c_2$, involved in respective estimates. In the theorem below  we provide the exact estimate for the ratio  (\ref{c1}).

\begin{theorem}\label{td4}
Assume that the L\'evy measure $\mu$  satisfies ($T_1$), ($T_2$),
and ($N_3$).

Then for every bounded set $A\subset \Re$
\begin{equation}\label{c3}
r_a(x)\sim e^{-a\zeta(x)},\quad x\to +\infty,
\end{equation}
uniformly in  $a\in A$.

In particular,
for any constants $c_1>c_*$ and $c_2<c_*$  there exists $y=y(c_1,c_2, A)$,  such that for  $x>y$, $a\in A$,
\begin{equation}\label{c311}
 x^{-c_1 a}\leq r_a(x)\leq x^{-c_2 a}
\ee
when $\mu$ is truncated, and
\begin{equation}\label{c32} x^{-c_1 a\ln^{-{1\over \beta}}x}\leq r_a(x)\leq x^{-c_2 a \ln^{-{1\over \beta}}x}
\ee
when $\mu$ is exponentially damped. Here $c_*$ is defined respectively in (\ref{cstar}) or (\ref{cstar1}).
\end{theorem}

\begin{proof}
By the inverse function theorem,
\be\label{der}
{d\over dx}\zeta(x)=\left[{d\over d\zeta}
\Mf_{1}(\zeta)\right]^{-1}\Big|_{\zeta=\zeta(x)}=
[\Mf_{2}(\zeta(x)]^{-1}.
\ee
Then
$$
{d\over dx}\ln \Kf (x)={\Mf_3(\zeta(x))\over \Mf_2(\zeta(x))}{d\over dx}\zeta(x)={\Mf_3(\zeta(x))\over \Mf_2^2(\zeta(x))}.
$$
If $\mu$ is supported in a bounded set, then $\Mf_3(\zeta)\sim (1/\sigma_+)\Mf_4(\zeta)$, $\zeta\to +\infty$ (see (\ref{731})). If $\mu$ is not supported in a bounded set, then $\Mf_3(\zeta)\ll \Mf_4(\zeta)$, $\zeta\to +\infty$ (see (\ref{733})). In both   cases, we have
$$
{d\over dx}\ln \Kf (x)\to 0, \quad x\to \infty
$$
because ($T_1$) provides ($\hat H_1$) (Lemma \ref{l22}). Thus,  for any bounded set $A$,
$$
{\Kf(x+a)\over \Kf(x)}=\exp\left(\int_x^{x+a} {d\over dy}\ln \Kf (y)\, dy\right)\to 0, \quad x\to
+\infty
$$
uniformly in $a\in A$. Therefore, by (\ref{answ6}),
$$
r_a(x)\sim e^{\Df(x+a)-\Df(x)},\quad x\to +\infty
$$
uniformly in  $a\in A$.

We have
$$
\Df(x+a)-\Df(x)=-a\zeta(x)-(x+a)[\zeta(x+a)-\zeta(x)]+\Mf_{0}(\zeta(x+a))-\Mf_{0}(\zeta(x)).
$$
Since  ${d\over d\zeta}\Mf_0(\zeta)=\Mf_1(\zeta)$ and $\Mf_1(\zeta(y))=y$, we get
$$
\ba
(x+a)&[\zeta(x+a)-\zeta(x)]-\Mf_{0}(\zeta(x+a))-\Mf_{0}(\zeta(x))=(x+a)\int_x^{x+a}\zeta'(y)\, dy-\int_x^{x+a}y\zeta'(y)\, dy\\&=\int_x^{x+a}(x+a-y)\zeta'(y)\, dy=\int_x^{x+a}\int_{y}^{x+a}\zeta'(y)\, dvdy=\int_0^{a}\int_{r}^{a}\zeta'(x+r)\, dsdr.
\ea
$$
Since  $\mu$ satisfies ($N_3$) we have (\ref{mu0}),  and therefore $\Mf_2(\zeta)\to +\infty$, $\zeta\to +\infty$. By (\ref{der}), this yields
$$
\zeta'(x)\to 0, \quad x\to +\infty.
$$
From the above relations we deduce that
$$
\Df(x+a)-\Df(x)\to -a\zeta(x), ,\quad x\to +\infty
$$
uniformly in  $a\in A$, which completes the proof of (\ref{c3}).

From (\ref{c3}) and (\ref{zeta0}) we deduce (\ref{c311}). From (\ref{c3}) and the first relation in (\ref{54}) we obtain (\ref{c32}).
\end{proof}

\textbf{Acknowledgement.}   The authors express their deep gratitude to the referee for inspiring  suggestions and   valuable remarks. We thank Dr. Kristian Evans, Swansea University, who made the proofreading of English in our paper.


\begin{thebibliography}{100}
\bibitem[Ku06]{Ku06} A. M. Kulik. Stochastic calculus of variations for general L\'evy processes
and its applications to jump-type SDE's  with non-degenerated
drift.  Available at   arxiv.org:math.PR/0606427v2.
\end{thebibliography}
\end{document}